%% file: semi_coherent_interfaces24-11__2_.tex
\documentclass[11pt,reqno]{amsart}
\usepackage[utf8]{inputenc}
\usepackage{amsmath,amssymb,amsthm,mathtools}
\usepackage[margin=1in]{geometry}
\usepackage{hyperref}
\usepackage{enumitem}
\usepackage[T1]{fontenc}
\usepackage{etoolbox}
\usepackage{xcolor}
\newcommand{\Chat}{\widehat{\mathcal C}}
\usepackage{tikz}
\usepackage{cite}

\numberwithin{equation}{section}

\makeatletter
\renewcommand{\@seccntformat}[1]{%
  \ifstrequal{#1}{subsection}
    {\textbf{\csname the#1\endcsname}\;\;}
    {\csname the#1\endcsname\quad}%
}
\makeatother

\theoremstyle{plain}

\newtheorem{theorem}[subsection]{Theorem}
\newtheorem{proposition}[subsection]{Proposition}
\newtheorem{lemma}[subsection]{Lemma}
\newtheorem{corollary}[subsection]{Corollary}

\theoremstyle{definition}
\newtheorem{definition}[subsection]{Definition}

\theoremstyle{remark}
\newtheorem{remark}[subsection]{Remark}

\newcommand{\C}{\mathcal{C}}
\newcommand{\eps}{\varepsilon}
\renewcommand{\d}{\textup{d}}
\newcommand{\E}{\mathcal{E}}

\newcommand{\spt}{\textup{spt}}
\newcommand{\Id}{\mathrm{Id}}
\newcommand{\F}{\mathcal{F}}
\newcommand{\Lip}{\operatorname{Lip}}
\newcommand{\ch}{\operatorname{ch}}

\newcommand{\Hu}{\mathcal{H}^1}
\newcommand{\Le}{\mathcal{L}^2}
\newcommand{\A}{\mathcal{A}}
\newcommand{\G}{\mathcal{G}}
\newcommand{\ov}{\overline}
\newcommand{\rb}{\bar\rho}
\newcommand{\res}{\mathbin{\vrule height 1.6ex depth 0pt width 0.13ex\vrule height 0.13ex depth 0pt width 1.3ex}}

\title[Energy of Semi-Coherent Interfaces]{On the Energy of Semi-Coherent Intercrystalline Interfaces }
\author{Valerio Corica}
\address{Sapienza Università di Roma}
\email{valerio.corica@uniroma1.it}
\email{emanuele.spadaro@uniroma1.it}

\author{Emanuele Spadaro}

\input{comandivalerio}

\date{\today}

\begin{document}

\begin{abstract}
	We develop an approach, in the spirit of Lauteri and Luckhaus (arXiv:1608.06155), to estimating the energy of semi-coherent interfaces between two crystals with different lattice spacings. Unlike classical approaches in physics and engineering, in this framework neither the geometry nor the location of the defects is prescribed \emph{a priori}; instead, both are optimized as part of the variational problem. Within this general geometric setting we prove matching upper and lower bounds of order $\beta |\log \beta|$ in the lattice misfit $\beta$, thus recovering the van der Merwe law for semi-coherent interfaces, the analog of the Read–Shockley scaling for grain boundaries.
    
\end{abstract}


\maketitle

\section{Introduction}

Semi-coherent interfaces arise when two crystals with the same lattice
structure but different lattice spacings meet, as for epitaxial thin films
grown on a crystalline substrate. For small mismatch, the misfit is
accommodated by a network of dislocations separating regions where the
two lattices are nearly coherent, and the energy of the interface results
from the competition between elastic relaxation and the cost of the defects.
If $\beta\in(0,1)$ denotes the lattice mismatch, classical continuum theories
predict, for $\beta\ll 1$, an interfacial energy
\begin{equation}\label{e:RS}
	\gamma \sim \beta|\log\beta|.
\end{equation}
This law was derived by Frank and van der Merwe for misfitting layers, and is analogous to the Read--Shockley formula for low-angle tilt grain boundaries~\cite{FrankVanDerMerwe1949a, FrankVanDerMerwe1949b, FrankVanDerMerwe1949c, ReadShockley1950, VanDerMerwe1950}. These computations, however, rely on strong assumptions on the geometry of the defects, such as a prescribed periodic arrangement of the dislocations.

In this paper we prove that \eqref{e:RS} holds in a variational model in
which neither the interface nor the defects are prescribed: their location,
their geometry and their possible periodicity are left to the minimization.
The model is a two-well extension of the geometric framework introduced by
Lauteri and Luckhaus~\cite{LauteriLuckhaus2016} for low-angle tilt grain
boundaries. Our main result, Theorem~\ref{t:main}, shows that the van der
Merwe scaling holds uniformly over all admissible interface and defect
geometries; in particular, the logarithmic factor is not an artifact of the
periodic ansatz.

\subsection{Background and motivation}

Semi-coherent interfaces have been studied through three main approaches.
Atomistic and Peierls--Nabarro-type models, going back to Frenkel and
Kontorova~\cite{FrenkelKontorova1938}, Peierls~\cite{Peierls1940} and
Nabarro~\cite{Nabarro1947} (see
\cite{DeLucaPonsiglioneSpadaro2024,FanzonPalombaroPonsiglione2017,FanzonPonsiglioneScala2020,FonsecaFuscoLeoniMorini2018,GinsterNeubauerZwicknagl2025,MullerPalombaro2013}
for recent mathematical developments), describe individual dislocations and,
in suitable regimes, their interactions and periodic arrangements. Continuum
elasticity models~\cite{FrankVanDerMerwe1949a,FrankVanDerMerwe1949b,FrankVanDerMerwe1949c,ReadShockley1950,VanDerMerwe1950}
treat the crystals as elastic media with dislocations in a prescribed
periodic configuration; they yield \eqref{e:RS} and show that, on
sufficiently large interfaces, a semi-coherent configuration is energetically
preferable to a fully coherent one. Finally, the geometric variational
approach of Lauteri and Luckhaus~\cite{LauteriLuckhaus2016} (see
also~\cite{LuckhausMugnai2010,LuckhausWohlgemuth2014} for its derivation from
mesoscopic many-body models) minimizes over both the lattice distortion and
the defect set, whose geometry is not fixed in advance, as expected in real
materials (see also \cite{ContiCrismaleGarroniMalusa2026,FortunaGarroniSpadaro2025, ScardiaTolotti2026}).

From the viewpoint of the calculus of variations, the periodicity of the
dislocations, their concentration on lower-dimensional sets and the formation
of interfaces should be consequences of energy minimization rather than
modeling assumptions.
For instance, periodicity was proved in~\cite{DeLucaPonsiglioneSpadaro2024} for a simplified Peierls--Nabarro model, and in~\cite{CoricaSantilli2026} for tilt grain boundaries. Here we take a complementary point of view: we derive sharp global energy bounds without assuming any specific
arrangement of the defects.

\subsection{The geometric two-well model}

We extend the Lauteri--Luckhaus framework to two wells, retaining the
geometric features of semi-coherent interfaces while avoiding a detailed
description of the atomic interactions.

Let $\varepsilon>0$ be the average lattice spacing. A configuration is a pair
$(A,S)$, where $S\subset\mathbb{R}^2$ is a closed set containing the defects
and $A\in L^1_{loc}(\mathbb{R}^2;\mathbb{R}^{2\times2})$ is the local
deformation of the lattice, whose incompatibility is concentrated on $S$:
\[
\operatorname{spt}(\curl A)\subset S,
\]
with $\curl A$ understood row-wise in the sense of distributions. Since
Burgers vectors are quantized in units of the lattice spacing, we require
that, for some $\tau_0>0$,
\begin{equation}\label{e:quantization}
\int_0^1 A(\gamma(t))\dot\gamma(t)\,\d t = 0 \qquad \textup{or} \qquad
\left|\int_0^1 A(\gamma(t))\dot\gamma(t)\,\d t\right| \geq \tau_0\,\varepsilon
\end{equation}
for every rectifiable closed curve $\gamma:[0,1]\to \R^2\setminus S$, the
circulation being understood in the weak sense.

The two crystals have lattice spacings $\varepsilon$ and $\varepsilon/\mu$,
with $\mu\in(0,1)$. Since $A$ maps the actual sample to a reference lattice,
their stress-free states are rotations and rotations composed with a
$\mu$-dilation, that is, the two wells are
\[
K_\mu:=\SO(2)\cup\mu\SO(2),
\]
and the mismatch parameter is $\beta:=1-\mu$. 
The energy of a configuration in a bounded open set
$\Omega\subset\mathbb{R}^2$ is the sum of the elastic energy and of the area
of the defect cores:
\[
\mathcal{F}_{\varepsilon}(A,S;\Omega)
:=
\int_{\Omega\setminus B_{\lambda\varepsilon}(S)}
\operatorname{dist}^2(A(x),K_\mu)\,\d x
+
\bigl|B_{\lambda\varepsilon}(S)\cap\Omega\bigr|,
\]
where $B_{\lambda\varepsilon}(S)$ is the $\lambda\varepsilon$-neighborhood of
$S$ and $\lambda>0$ is the core radius in units of the lattice spacing. The
energy is rotation invariant and does not prescribe the geometry of $S$: its
purpose is not to describe the microscopic structure in detail, but to
isolate the mechanisms responsible for the global cost of accommodating two
incompatible lattices.

\subsection{Main result}

Fix $\delta>\varepsilon>0$. We denote by
$\mathcal{C}_\varepsilon(\mu)$ the class of configurations satisfying
\eqref{e:quantization} and
\begin{equation}\label{e:prescribed-phase}
A(x)=
\begin{cases}
	\Id,
	& x\in\R\times(\delta,+\infty),\\[1mm]
	\mu\Id,
	& x\in\R\times(-\infty,-\delta),
\end{cases}
\end{equation}
so that the two crystals occupy the half-planes above and below the strip
$\R\times(-\delta,\delta)$. Every admissible configuration must interpolate
between the two wells, but the location and the geometry of the transition
region, as well as the defect set, are free. For $L>\delta$ we set
$Q_L:=(-L,L)^2$ and
\begin{equation}\label{e:E_L}
	E_L:=
	\liminf_{\varepsilon\to0^+}
	\inf_{(A,S)\in\mathcal{C}_\varepsilon(\mu)}
	\eps^{-1}
	\mathcal{F}_{\varepsilon}(A,S;Q_L).
\end{equation}

\begin{theorem}\label{t:main}
	There exist $\mu_0\in(0,1)$ and $C>0$, depending on $\tau_0$ and $\lambda$,
	such that for every $\mu\in(\mu_0,1)$, $\delta>0$ and $L>\delta$
	\[
	C^{-1}L\beta|\log\beta|
	\leq E_L
	\leq CL\beta|\log\beta|,
	\qquad \beta=1-\mu.
	\]
\end{theorem}

Theorem~\ref{t:main} is the analogue for semi-coherent interfaces of the
Read--Shockley lower bound proved in~\cite{LauteriLuckhaus2016} for low-angle
tilt grain boundaries. There, the field takes values in the single well
$\SO(2)$ and the interface separates two grains with different orientations. Here the two phases lie
in different wells and the interface is forced by \eqref{e:prescribed-phase}
rather than by a rotation mismatch, which calls for new ingredients at each
step of the proof (see Section~\ref{ss:ideas}). The upper bound is obtained
by a periodic arrangement of defects along a line, at mutual distance of
order $\eps/\beta$ (Appendix~\ref{app:upper-bound}; see
also~\cite{FortunaGarroniSpadaro2025}).

\subsection{Ideas of the proof}\label{ss:ideas}

The main difficulty is that a field with small elastic energy may be close to
either $\SO(2)$ or $\mu\SO(2)$: the lower bound must control at the same time
the choice of the well, the transition between the two wells and the
incompatibility generated by an arbitrary defect set. We follow the overall
strategy of Lauteri and Luckhaus, with the following new ingredients.

In Section~\ref{s:uniform} we prove an interior rigidity estimate for
gradients close to $K_\mu$, with a constant uniform in $\mu\in[1/2,2]$. Since
the two wells merge as $\beta\to0$, this uniformity is what makes the
constants of Theorem~\ref{t:main} independent of the misfit.

In Section~\ref{s:rigidity} we prove a rigidity estimate for incompatible
fields with multiple wells (Theorem~\ref{t:rigidity}), valid in any dimension
for strongly incompatible wells, by extending the blow-up argument
of~\cite{DeLellisSzekelyhidi2006}. Its constant may
depend on the wells, which is harmless since it is only used for compactness
at fixed $\mu$.

In Section~\ref{s:transfer} we prove the energy transfer estimate; here the
harmonic replacement relies on the conformal structure of $K_\mu$ (see the
discussion preceding Lemma~\ref{l:harmonic}).

In Section~\ref{s:log} we isolate an abstract version of the logarithmic
estimate of~\cite{LauteriLuckhaus2016} (Theorem~\ref{t:log}), in which the
only input from the elastic model is the energy transfer inequality; we
believe this formulation may be useful in other settings. 

The lower bound is finally proved in Section~\ref{s:lower}. The logarithmic enhancement arises from the competition, across all relevant length scales, between the elastic energy induced by the defects and the cost of their cores.

\section{Uniform rigidity for conformal double wells}\label{s:uniform}

In order to obtain constants in Theorem~\ref{t:main} which do not degenerate in the small-misfit regime, we need a rigidity estimate as Friesecke--James--M\"uller \cite{Friesecke-James-Muller2002} for $K_\mu$ whose constant is independent of $\mu$ in a neighborhood of $\mu=1$.
The rigidity estimates in \cite{ChaudhuriMuller2003,DeLellisSzekelyhidi2006, FaracoZhong2005, LamyLorentPeng2023} involve constants depending on the separation of the two wells and do not provide the uniformity required here.
The main result of this section is the following estimate.

\begin{theorem}[Uniform rigidity estimate]\label{t:uniform rigidity}
Let $\Omega\subset\R^2$ be a bounded, connected Lipschitz domain. Then there exists a constant $C_U(\Omega)>0$ with the following property.
	Let $\mu\in[\frac12,2]$, and let $u\in W^{1,2}(\Omega;\R^2)$. Then there exists $S\in K_\mu$ such that
	\[
	\int_{\Omega}|\nabla u-S|^2\,\d x\le C_U\int_{\Omega}\dist_{K_\mu}^2(\nabla u)\,\d x .
	\]
\end{theorem}

The proof follows the scheme of the classical rigidity estimate of Friesecke, James and M\"uller \cite{Friesecke-James-Muller2002}, as adapted to two wells by Chaudhuri and M\"uller \cite{ChaudhuriMuller2003}.
The key for the uniformity in $\mu$ varying in a compact set of $(0,+\infty)$ is a linear rigidity estimate for scalar harmonic functions close to two values (Lemma~\ref{l:h-rigidity}), which is applied to $\log |F|$, where $F$ is the projection of the gradient on the conformal matrices.
In particular, the case $\mu = 1$ is included: we thus recover an alternative proof of the classical two-dimensional result, which does not use Korn’s inequality.

\subsection{Algebra of conformal matrices}
We use $v\cdot w:=\sum_i v_iw_i$ for the Euclidean scalar product
and $A:B:=\sum_{i,j}A_{ij}B_{ij}$ for the Frobenius scalar product,
with the corresponding norms denoted by $|\cdot|$.
We denote by $\Id$ the identity matrix and by $J$ the counterclockwise rotation by $\frac\pi2$:
\[
\Id = \begin{pmatrix}
	1 & 0 \\
	0 & 1
\end{pmatrix}, \qquad
J = \begin{pmatrix}
	0 & -1 \\
	1 & 0
\end{pmatrix}.
\]
The conformal matrices (including $0$) form the two-dimensional linear subspace
\[
\operatorname{CO}_+ := \left\{a \Id + b J : a, b\in \R \right\}
=\left\{\begin{pmatrix}
	a & -b \\
	b & a
\end{pmatrix}: a,b\in\R\right\}.
\]
For $A = \begin{pmatrix}
	a_{11} & a_{12} \\
	a_{21} & a_{22}
\end{pmatrix}$ the cofactor matrix is
\[
\operatorname{cof}(A) = \begin{pmatrix}
	a_{22} & -a_{21} \\
	-a_{12} & a_{11}
\end{pmatrix}.
\]
In dimension two the map $A\mapsto\operatorname{cof}(A)$ is linear, isometric for the Frobenius norm, involutive and self-adjoint, i.e.\ $\operatorname{cof}(A):B=A:\operatorname{cof}(B)$.
Moreover,
\begin{equation}
    \operatorname{cof}(A)= A 
\end{equation} if and only if $A \in\operatorname{CO}_+$. 
This is specific to dimension two: in higher dimensions the equation is satisfied only by matrices in $\operatorname{SO}(n)$ (and by $A = 0$).
We set
\[
PA := \frac{A+\operatorname{cof}(A)}{2} \qquad \text{and}\qquad P^\perp A := \frac{A-\operatorname{cof}(A)}{2}
\]
for the orthogonal projections onto $\operatorname{CO}_+$ and on its orthogonal complement (the anti-conformal matrices), respectively.
Since $K_\mu\subset \operatorname{CO}_+$, we get
\begin{equation}\label{e:projection}
	|A-\operatorname{cof}(A)| = 2|P^\perp A| = 2 \dist_{\operatorname{CO}_+}(A) \leq 2 \dist_{K_\mu}(A)\qquad\forall\,A\in\R^{2\times2}.
\end{equation}

For $F=a\Id+bJ\in\operatorname{CO}_+$ we set
\[
\varrho(F):=\sqrt{a^2+b^2}=\frac{|F|}{\sqrt2}.
\]
For $\lambda>0$ and $\theta\in\R$ we have
\[
\big|F-\lambda(\cos\theta\, \Id+\sin\theta\, J)\big|^2=2\big(\varrho(F)^2+\lambda^2-2\lambda(a\cos\theta+b\sin\theta)\big)\ge2\big(\varrho(F)-\lambda\big)^2,
\]
with equality for a suitable $\theta$. Hence
\begin{equation}\label{e:distK}
	\dist_{\lambda\SO(2)}(F)=\sqrt2\,|\varrho(F)-\lambda|,
	\qquad
	\dist_{K_\mu}(F)=\sqrt2\,\dist\big(\varrho(F),\{1,\mu\}\big)\qquad\forall\,F\in\operatorname{CO}_+ .
\end{equation}
Finally, since $|\mu R-\nu R|=\sqrt2\,|\mu-\nu|$ for every $R\in\SO(2)$, the Hausdorff distance between $K_\mu$ and $K_\nu$ is at most $\sqrt2|\mu-\nu|$, and therefore
\begin{equation}\label{e:dist-mu}
	\big|\dist_{K_\mu}(A)-\dist_{K_\nu}(B)\big|\le|A-B|+\sqrt2\,|\mu-\nu|\qquad\forall\,A,B\in\R^{2\times2},\ \mu,\nu>0 .
\end{equation}

\subsection{Harmonic conformal fields}
Divergence and curl of matrix fields are taken row-wise.
For a conformal field $F=a\Id+bJ$, with rows $(a,-b)$ and $(b,a)$, we have
\[
\operatorname{div}F=\begin{pmatrix}\partial_1a-\partial_2b\\ \partial_1b+\partial_2a\end{pmatrix},
\qquad
\curl F=\begin{pmatrix}-\partial_1b-\partial_2a\\ \partial_1a-\partial_2b\end{pmatrix}.
\]
In particular, for conformal fields the conditions $\operatorname{div}F=0$ and $\curl F=0$ are equivalent, and both are equivalent to the Cauchy--Riemann equations
\begin{equation}\label{e:CR}
	\partial_1a=\partial_2b,\qquad \partial_2a=-\partial_1b .
\end{equation}

\begin{definition}\label{d:harmonic conformal}
	Let $U\subset\R^2$ be open. A \emph{harmonic conformal field} in $U$ is a field $F=a\Id+bJ\in L^1_{\rm loc}(U;\operatorname{CO}_+)$ such that $\operatorname{div}F=0$ (equivalently, $\curl F=0$) in the sense of distributions in $U$.
\end{definition}
\begin{remark}
By \eqref{e:CR}, $\Delta a=\Delta b=0$ in distributions, so that $a$ and $b$ are harmonic functions, hence smooth.
Moreover \eqref{e:CR} gives $\nabla b=J\nabla a$, i.e.\ $|\nabla a|=|\nabla b|$ and $\nabla a\cdot\nabla b=0$. Setting $s:=a^2+b^2=\varrho(F)^2$, we deduce the identities
\begin{equation}\label{e:identities}
	|DF|^2=2|\nabla a|^2+2|\nabla b|^2=4|\nabla a|^2,\qquad
	\Delta s=|DF|^2,\qquad
	|\nabla s|^2=4s|\nabla a|^2=s\,|DF|^2 ,
\end{equation}
where we used $\Delta s=2|\nabla a|^2+2|\nabla b|^2+2a\Delta a+2b\Delta b$ and $\nabla s=2a\nabla a+2b\nabla b$.
On the open set $\{s>0\}$ we define $h:=\frac12\log s=\log\varrho(F)$. By \eqref{e:identities},
\begin{equation}\label{e:log-harmonic}
	\Delta h=\frac{\Delta s}{2s}-\frac{|\nabla s|^2}{2s^2}=\frac{|DF|^2}{2s}-\frac{|DF|^2}{2s}=0,
	\qquad
	|DF|^2=4s\,|\nabla h|^2=4\varrho(F)^2|\nabla h|^2 .
\end{equation}
\end{remark}
The next lemma provides the harmonic conformal competitor.

\begin{lemma}\label{l:harmonic_conformal}
	Let $\Omega\subset \R^2$ be a bounded open set and $\mu>0$.
	For $u\in W^{1,2}(\Omega;\R^2)$, there exist $w\in u + H^1_0(\Omega;\R^2)$ harmonic, and $F\in L^2(\Omega;\operatorname{CO}_+)$ a harmonic conformal field such that
	\begin{gather}
		\|\nabla u - \nabla w\|_{L^2(\Omega)}\leq C\|\dist_{K_\mu} (\nabla u)\|_{L^2(\Omega)},\label{e:h_c-1}\\
		\|\nabla u - F\|_{L^2(\Omega)}\leq C\|\dist_{K_\mu} (\nabla u)\|_{L^2(\Omega)},\label{e:harm_conf-1}\\
		\|\dist_{K_\mu} (F)\|_{L^2(\Omega)} \leq C\|\dist_{K_\mu} (\nabla u)\|_{L^2(\Omega)},\label{e:harm_conf-2}
	\end{gather}
	with $C>0$ a universal constant.
\end{lemma}

\begin{proof}
By Piola's identity, for every $v\in W^{1,2}(\Omega;\R^2)$
\begin{equation}\label{e:piola}
		\int_\Omega\operatorname{cof}(\nabla v):\nabla\varphi\,\d x=0\qquad\forall\,\varphi\in H^1_0(\Omega;\R^2).
	\end{equation}
Let $v$ be the solution to
\begin{equation}
    \begin{cases}
        \Delta v = \operatorname{div} ( \operatorname{cof}( \nabla u)-\nabla u  ) \quad \text{ in } \Omega, \\
        
        v \in H^1_0(\Omega;\R^2).
    \end{cases}
\end{equation}
Then $w=u+v$ is the harmonic replacement, and testing the equation with $v$ gives
	\begin{align*}
	\|\nabla u- \nabla w\|_{L^2(\Omega)}=	\|\nabla v\|_{L^2(\Omega)} \leq \|\operatorname{cof}( \nabla u)-\nabla u  \|_ {L^2(\Omega)} \leq 2\|\dist_{K_\mu} (\nabla u)\|_{L^2(\Omega)}
	\end{align*}
	which gives \eqref{e:h_c-1}.
	Set $F:=P\nabla w$. By Piola's identity and the harmonicity of $w$,
    $\operatorname{div}F=\frac12\Delta w+\frac12\operatorname{div}(\operatorname{cof}\nabla w)=0$,
    so $F$ is a harmonic conformal field.
    Next, $\nabla u-F=P^\perp\nabla u+P(\nabla u-\nabla w)$, and $P$ is an orthogonal projection; hence, by \eqref{e:projection} and \eqref{e:h_c-1},
	\[
	\|\nabla u - F\|_{L^2(\Omega)} \leq \|P^\perp \nabla u\|_{L^2(\Omega)} + \|\nabla u - \nabla w\|_{L^2(\Omega)}\le 3 \|\dist_{K_\mu} (\nabla u)\|_{L^2(\Omega)}.
	\]
	Finally, \eqref{e:harm_conf-2} follows from \eqref{e:harm_conf-1} and the triangle inequality.
\end{proof}

\subsection{Rigidity of harmonic conformal fields}

\begin{lemma}\label{l:h_c-rigidity-1}
	Let $U\subset \R^2$ be a connected open set, $\mu>0$, and let $F$ be a harmonic conformal field in $U$ with $F(x)\in K_\mu$ for a.e.\ $x\in U$. Then $F$ is constant.
\end{lemma}

\begin{proof}
	Since $F$ is continuous and $K_\mu$ is closed, $F(x)\in K_\mu$ for every $x\in U$. Hence the continuous function $s=\varrho(F)^2$ takes values in $\{1,\mu^2\}$ and, $U$ being connected, it is constant. By \eqref{e:identities}, $|DF|^2=\Delta s=0$, and $F$ is constant.
\end{proof}

Note that the lemma includes the case $\mu=1$, and that no quantitative rigidity estimate is used in its proof. Its quantitative counterpart is the following.

\begin{lemma}\label{l:h_c-rigidity-2}
	There exists a universal constant $\eps_0>0$ with the following property.
	Let $\mu\in [\frac12,2]$ and let $F$ be a harmonic conformal field in $B_2$ with
	\[
	\|\dist_{K_\mu}(F)\|_{L^2(B_2)}\leq \eps_0 .
	\]
	Then $\frac14\le\varrho(F(x))\le4$ for every $x\in \overline{B}_{7/4}$.
\end{lemma}

\begin{proof}
	We argue by contradiction, and assume that there exist $\mu_j\in[\frac12,2]$, harmonic conformal fields $F_j$ in $B_2$ and points $x_j\in\overline{B}_{7/4}$ such that
	\[
	\|\dist_{K_{\mu_j}}(F_j)\|_{L^2(B_2)}\to 0 \qquad\text{and}\qquad \varrho(F_j(x_j))\notin\big[\tfrac14,4\big].
	\]
	Since $|F|\le\dist_{K_\mu}(F)+2\sqrt2$ for $\mu\le2$, the sequence $(F_j)$ is bounded in $L^2(B_2)$.
	By the interior estimates for harmonic functions, up to a subsequence $F_j\to F$ in $C^1_{\rm loc}(B_2)$, $\mu_j\to\mu\in[\frac12,2]$ and $x_j\to\bar x\in\overline{B}_{7/4}$. The limit $F$ is conformal and satisfies \eqref{e:CR}, i.e.\ it is a harmonic conformal field in $B_2$.
	By \eqref{e:dist-mu}, $\dist_{K_{\mu_j}}(F_j)\to\dist_{K_\mu}(F)$ locally uniformly in $B_2$, so that
	\[
	\int_{B_s}\dist^2_{K_\mu}(F)\,\d x=\lim_{j}\int_{B_s}\dist^2_{K_{\mu_j}}(F_j)\,\d x=0\qquad\forall\,s<2 .
	\]
	By Lemma~\ref{l:h_c-rigidity-1}, $F$ is constant with $\varrho(F)\in\{1,\mu\}\subset[\frac12,2]$.
	On the other hand, the uniform convergence on $\overline{B}_{7/4}$ gives $\varrho(F_j(x_j))\to\varrho(F(\bar x))$, and therefore $\varrho(F(\bar x))\in[0,\frac14]\cup[4,\infty)$, a contradiction.
\end{proof}

\begin{remark}\label{r:no separation}
	The threshold $\eps_0$ is uniform also as the two wells merge: the limit $\mu$ in the proof may be equal to $1$, and Lemma~\ref{l:h_c-rigidity-1} does not require $\mu\neq1$. Only the bounds $\frac12\le\mu\le2$ are used.
\end{remark}

\subsection{Rigidity for scalar harmonic functions close to two values}

\begin{lemma}\label{l:h-rigidity}
	For every $0<r<R$ there exists $C(r,R)>0$ with the following property: for every $a\in \R$ and for every harmonic function $h:B_R\to \R$,
	\begin{equation}\label{e:h-rigidity}
		\min_{c\in \{0,a\}} \|h-c\|_{L^2(B_r)}\leq C(r,R)\, \|\dist(h, \{0,a\})\|_{L^2(B_R)}.
	\end{equation}
	The constant does not depend on $a$; in particular, it does not degenerate as $a\to0$.
\end{lemma}

\begin{proof}
	If $a=0$, then $\dist(h,\{0,a\})=|h|$ and \eqref{e:h-rigidity} holds with $C=1$.
	If $a\ne0$, the function $g:=h/a$ is harmonic; since both sides of \eqref{e:h-rigidity} are scaling invariant, it suffices to consider $a=1$.

	Assume by contradiction that for every $j\in\N$ there exists a harmonic function $h_j:B_R\to\R$ with
	\begin{equation}\label{e:contra-h}
		\min_{c\in \{0,1\}} \|h_j-c\|_{L^2(B_r)} > j\, \|\dist(h_j, \{0,1\})\|_{L^2(B_R)}.
	\end{equation}
	Since $|h_j|\le\dist(h_j,\{0,1\})+1$, the left-hand side of \eqref{e:contra-h} is at most $\|\dist(h_j, \{0,1\})\|_{L^2(B_R)}+|B_R|^{1/2}$. Hence
	\[
	(j-1)\,\|\dist(h_j, \{0,1\})\|_{L^2(B_R)}<|B_R|^{1/2},
	\]
	so that $\|\dist(h_j, \{0,1\})\|_{L^2(B_R)}\to0$ and $(h_j)$ is bounded in $L^2(B_R)$.
	By the interior estimates for harmonic functions, up to a subsequence $h_j\to h$ locally uniformly in $B_R$, with $h$ harmonic and $\dist(h,\{0,1\})=0$ in $B_R$. Since $B_R$ is connected and $h$ is continuous, $h\equiv c$ for some $c\in\{0,1\}$.
	By the uniform convergence on $\overline B_r$, for $j$ large we have $|h_j-c|<\frac12$ on $B_r$, hence $|h_j-c|=\dist(h_j,\{0,1\})$ on $B_r$, and
	\[
	\min_{c'\in \{0,1\}} \|h_j-c'\|_{L^2(B_r)}\le\|h_j-c\|_{L^2(B_r)}\le\|\dist(h_j, \{0,1\})\|_{L^2(B_R)},
	\]
	which contradicts \eqref{e:contra-h} for $j\ge1$.
\end{proof}

\subsection{Proof of Theorem~\ref{t:uniform rigidity}}\label{ss:proof uniform}
In Steps 1 and 2, $C$ denotes a universal constant which may change from line to line.

\medskip

\noindent\emph{Step 1: interior estimate for harmonic conformal fields.}
We show that for every $\mu\in[\frac12,2]$ and every harmonic conformal field $F\in L^2(B_2;\operatorname{CO}_+)$
\begin{equation}\label{e:rigidity_F}
	\|DF\|_{L^2(B_1)} \leq C\, \|\dist_{K_\mu}(F)\|_{L^2(B_2)}.
\end{equation}
Set $e:=\|\dist_{K_\mu}(F)\|_{L^2(B_2)}$ and let $\eps_0$ be the constant of Lemma~\ref{l:h_c-rigidity-2}.

If $e>\eps_0$, we apply the Caccioppoli inequality to the harmonic components of $F$ and use $|F|\le\dist_{K_\mu}(F)+2\sqrt2$:
\[
\|DF\|_{L^2(B_1)} \leq C\, \|F\|_{L^2(B_2)} \leq C\big(e + 2\sqrt2\,|B_2|^{1/2}\big) \leq C \left(1+\frac{1}{\eps_0}\right)e.
\]

If $e\le\eps_0$, by Lemma~\ref{l:h_c-rigidity-2} we have $\frac14\le\varrho(F)\le4$ on $\overline{B}_{7/4}$. Hence the function $h:=\log\varrho(F)$ is well defined and harmonic in $B_{7/4}$ by \eqref{e:log-harmonic}.
Both $\varrho(F)$ and the values $\lambda\in\{1,\mu\}$ belong to $[\frac14,4]$, where $\log$ is $4$-Lipschitz; therefore, by \eqref{e:distK},
\[
\dist\big(h,\{0,\log\mu\}\big)\le4\,\dist\big(\varrho(F),\{1,\mu\}\big)=2\sqrt2\,\dist_{K_\mu}(F)\qquad\text{in }B_{7/4}.
\]
Lemma~\ref{l:h-rigidity}, applied with $r=\frac32$, $R=\frac74$ and $a=\log\mu$, gives $c\in\{0,\log\mu\}$ with
\[
\|h-c\|_{L^2(B_{3/2})}\le C\,\|\dist(h,\{0,\log\mu\})\|_{L^2(B_{7/4})}\le C e .
\]
Since $h-c$ is harmonic, the Caccioppoli inequality and \eqref{e:log-harmonic}, together with $\varrho(F)\le4$, give
\begin{equation}\label{e:balls}
\|DF\|_{L^2(B_1)}\le 8\,\|\nabla h\|_{L^2(B_1)}\le C\,\|h-c\|_{L^2(B_{3/2})}\le C e .
\end{equation}
This proves \eqref{e:rigidity_F} in both cases. Note that the case $\mu=1$ corresponds to $a=0$ in Lemma~\ref{l:h-rigidity}.

\medskip

\noindent\emph{Step 2: weighted estimate.}
Now the goal is to show that for any bounded, connected Lipschitz domain $\Omega\subset\R^2$ and $F \in L^2(\Omega;CO_+)$ harmonic and conformal,
\begin{equation}\label{e:weighted}
    \int_{\Omega} d_{\partial\Omega}^2 |DF(x)|^2 \,\d x \le C \int_{\Omega}\dist_{K_\mu}(F)^2 \,\d x.
\end{equation}

By the scaling and translation invariance the interior estimate \eqref{e:rigidity_F} can be extended to any pair of concentric balls $B_r\subset B_{2r}\subset\Omega$, and as a consequence one can obtain the same estimate for any couple of concentric cubes in the following form
\begin{equation}
	\|DF\|_{L^2(Q_L)} \leq \frac{C}{L}\, \|\dist_{K_\mu}(F)\|_{L^2(Q_{4L})}.
\end{equation}

By the Whitney decomposition theorem, there exists a countable
family of dyadic cubes $\{Q_j\}_j$ with pairwise disjoint
interiors such that, writing $\ell(Q_j)$ for the side length of $Q_j$,
\[
\Omega=\bigcup_j Q_j,
\qquad
4\sqrt{2}\,\ell(Q_j)
\leq \operatorname{dist}(Q_j,\partial\Omega)
\leq 10\sqrt{2}\,\ell(Q_j).
\]
Moreover, $4Q_j\subset\Omega$ and the dilated cubes satisfy
the uniform bounded overlap property
\[
\sum_j \mathbf{1}_{4Q_j}(x)\leq N
\qquad\text{for every }x\in\mathbb R^2.
\]

Therefore, we obtain 
\begin{equation*}
    \int_{Q_j} \frac{d^2_{\partial\Omega}}{\ell(Q_j)^2}|DF(x)|^2 \,\d x \leq \frac{C}{\ell(Q_j)^2} \int_{4Q_j}\dist_{K_\mu}(F)^2 \,\d x.
\end{equation*}
Multiplying by $\ell(Q_j)^2$ and summing on the cubes, by the bounded overlap property, we obtain \eqref{e:weighted}.

\medskip

\noindent\emph{Step 3: conclusion.}
Let $u\in W^{1,2}(\Omega;\R^2)$, set $E:=\|\dist_{K_\mu}(\nabla u)\|_{L^2(\Omega)}$, and let $F$ be the harmonic conformal field given by Lemma~\ref{l:harmonic_conformal} on $\Omega$. In this step $C$ may depend on $\Omega$. By the weighted Poincar\'e inequality on $\Omega$ (see \cite{Friesecke-James-Muller2002}), \eqref{e:weighted} and \eqref{e:harm_conf-2},
\[
\|F-(F)_{\Omega}\|_{L^2(\Omega)} \leq C\, \|d_{\partial\Omega}DF\|_{L^2(\Omega)} \leq C\, \|\dist_{K_\mu}(F)\|_{L^2(\Omega)}\leq C E.
\]
Since subtracting the mean value is an orthogonal projection in $L^2(\Omega)$, by \eqref{e:harm_conf-1} we get
\[
\|\nabla u-(\nabla u)_{\Omega}\|_{L^2(\Omega)} \leq \|\nabla u-F\|_{L^2(\Omega)} + \|F-(F)_{\Omega}\|_{L^2(\Omega)}\leq C E.
\]
Let $S\in K_\mu$ be a nearest point to $(\nabla u)_{\Omega}$ in the compact set $K_\mu$. Integrating the inequality $\dist_{K_\mu}((\nabla u)_{\Omega})\le|(\nabla u)_{\Omega}-\nabla u(x)|+\dist_{K_\mu}(\nabla u(x))$ over $\Omega$, we obtain
\[
|\Omega|^{1/2}\,|(\nabla u)_{\Omega} - S| \le \|\nabla u-(\nabla u)_{\Omega}\|_{L^2(\Omega)}+E\le CE .
\]
Therefore
\[
\|\nabla u-S\|^2_{L^2(\Omega)} = \|\nabla u-(\nabla u)_{\Omega}\|^2_{L^2(\Omega)} + |\Omega|\,|(\nabla u)_{\Omega} - S|^2\leq C E^2,
\]
which is the claim with $C_U=C_U(\Omega)$ independent of $\mu$.
Since the proof is invariant under rescaling, the constants depending on $\Omega$ depend only on its shape and are unchanged under uniform dilations. \qed

\begin{remark}\label{r:uniform comments}
	We collect some comments on Theorem~\ref{t:uniform rigidity}.
	\begin{itemize}
		\item[(i)] The same proof gives, for every $0<m\le M<\infty$, a constant $C_U(\Omega,m,M)$ valid for all $\mu\in[m,M]$: it suffices to replace the interval $[\frac14,4]$ in Lemma~\ref{l:h_c-rigidity-2} with $[\frac12\min\{1,m\},2\max\{1,M\}]$. For $\mu\in[\frac12,2]$ the constant depends only on $\Omega$. The same argument holds for a finite number of $\mu_i \in [m,M]$.
		\item[(ii)] For $\mu=1$ the theorem gives the two-dimensional rigidity estimate of \cite{Friesecke-James-Muller2002}, and the argument above gives a proof of it.
	\end{itemize}
\end{remark}

\section{Rigidity estimate for incompatible fields and multiple wells}\label{s:rigidity}

In dislocation theory, rigidity estimates of the type proved in the previous section must also account for the curl of the strain field. In two dimensions, Garroni, Leoni and Ponsiglione \cite{GarroniLeoniPonsiglione2010} proved a Korn inequality for incompatible fields, and M\"uller, Scardia and Zeppieri \cite{MullerScardiaZeppieri2014} established its nonlinear counterpart for the well $\SO(2)$. Conti and Garroni \cite{ContiGarroni2021} obtained sharp versions of both estimates in higher dimensions. In the linear setting, Gmeineder, Lewintan and Neff \cite{GmeinederLewintanNeff2023} characterized the linear maps for which a Korn-type estimate with curl control holds in arbitrary dimension.

In multiwell problems, 
Almi, Reggiani and Solombrino \cite{AlmiReggianiSolombrino2026} proved a rigidity estimate for incompatible fields near multiple wells that involves both the curl and the divergence of the field. For deformation gradients near two rank-one connected wells, Davoli and Friedrich \cite{DavoliFriedrich2020} obtained a rigidity estimate involving the second derivatives of the deformation.

Here we consider two strongly incompatible wells. We prove that the $L^{1^*}$-distance of a field from one of the two wells is controlled by its $L^{1^*}$-distance from their union and the total variation of its curl. This estimate yields the compactness needed to extend the Lauteri--Luckhaus theory to semi-coherent interfaces.

\subsection{Young measures and strongly incompatible wells}

Throughout this subsection $n\ge2$, $d,m\geq1$, $p\in [1,+\infty]$ and $\Omega\subset\R^n$ is a bounded open set.
We recall some definitions and facts about Young measures $\Y^p(\Omega;\R^d)$ and gradient Young measures $\GY^p(\Omega;\R^{m\times n})$; we refer to \cite{Muller1999} for the proofs.

\begin{itemize}
	\item[(Y1)] A sequence $(v_j)$ bounded in $L^p(\Omega;\R^d)$ \emph{generates} the Young measure $\nu=\{\nu_x\}_{x\in \Omega}$ if $\nu_x$ is a probability measure on $\R^d$ for a.e.\ $x\in\Omega$, the map $x\mapsto\nu_x$ is weakly$^*$ measurable, and
	\[
	\lim_{j\to +\infty} \int_\Omega f(x, v_j(x))\, \d x = \int_\Omega \int_{\R^d}f(x, P)\,  \d \nu_x(P) \, \d x
	\]
	for every Carath\'eodory function $f:\Omega\times \R^d\to \R$ such that $(f(\cdot, v_j))_j$ is equi-integrable. We write $\langle\nu_x,f\rangle:=\int_{\R^d}f(x,P)\,\d\nu_x(P)$, and $\nu\in\Y^p(\Omega;\R^d)$.

	\item[(Y2)] If $\nu$ is generated by a sequence $(\nabla u_j)$, with $(u_j)$ bounded in $W^{1,p}(\Omega;\R^m)$, then $\nu$ is a \emph{gradient Young measure}, and we write $\nu\in\GY^p(\Omega;\R^{m\times n})$.

	\item[(Y3)] \emph{Fundamental theorem.} Every bounded sequence in $L^p(\Omega;\R^d)$ has a subsequence which generates a Young measure. Moreover, two sequences $(v_j)$ and $(w_j)$ bounded in $L^p$ with $v_j-w_j\to0$ in measure generate the same Young measure.
\end{itemize}

\begin{definition}\label{d:incompatible}
	Two disjoint compact sets $K_1,K_2\subset\R^{m\times n}$ are \emph{strongly incompatible} if for every connected open set $\Omega'\subset\R^n$ and every $\nu\in\GY^\infty(\Omega';\R^{m\times n})$ with $\supp\nu_x\subset K_1\cup K_2$ for a.e.\ $x\in\Omega'$, either $\supp\nu_x \subset K_1$ for a.e.\ $x\in \Omega'$, or $\supp\nu_x \subset K_2$ for a.e.\ $x\in \Omega'$.
\end{definition}

For a compact set $K\subset\R^{m\times n}$ we write
\[
\dist_K(M):=\min_{P\in K}|M-P|.
\]
We will use the following elementary facts.

\begin{lemma}\label{l:pointwise}
	Let $K_1,K_2\subset\R^{m\times n}$ be disjoint compact sets, $K:=K_1\cup K_2$, and set
	\[
	\eta:=\dist(K_1,K_2)>0,\qquad D:=\max\Big\{\max_{P\in K_1}\dist_{K_2}(P),\,\max_{P\in K_2}\dist_{K_1}(P)\Big\},\qquad C_K:=\max_{P\in K}|P| .
	\]
	Then for every $M\in\R^{m\times n}$, every $q\ge1$ and $i\in\{1,2\}$:
	\begin{enumerate}
		\item[(a)] $|M|\le \dist_K(M)+C_K$;
		\item[(b)] $\dist_{K_i}(M)\leq \dist_K(M)+D$, and hence $\dist_{K_i}^{q}(M)\leq 2^{q-1}\big(\dist_K^{q}(M)+D^{q}\big)$;
		\item[(c)] $\dist_{K_1}(M)+\dist_{K_2}(M)\ge\eta$, and hence $\dist_{K_1}^{q}(M)+\dist_{K_2}^{q}(M)\ge2^{1-q}\eta^{q}$.
	\end{enumerate}
\end{lemma}

\begin{proof}
	(a) is immediate. For (b), let $P\in K$ with $|M-P|=\dist_K(M)$; then $$\dist_{K_i}(M)\le|M-P|+\dist_{K_i}(P)\le\dist_K(M)+D.$$ For (c), if $P_i\in K_i$ realize $\dist_{K_i}(M)$, then $\eta\le|P_1-P_2|\le \dist_{K_1}(M)+\dist_{K_2}(M)$. The inequalities for the powers follow from the convexity of $t\mapsto t^q$.
\end{proof}

\subsection{Rigidity estimate}

Throughout this subsection $n\ge2$, $m\ge1$, $\Omega\subset\R^n$ is a bounded open set and $1^*:=\frac{n}{n-1}$.
For $A\in L^1_{\rm loc}(\Omega;\R^{m\times n})$, $\curl A$ is the tensor
\[
(\curl A)_{ijk}=\partial_jA_{ik}-\partial_kA_{ij},\qquad i=1,\dots,m,\ j,k=1,\dots,n,
\]
with derivatives in the sense of distributions. We denote by $|\curl A|$ its total variation when $\curl A$ is a finite Radon measure.
For a compact set $K\subset\R^{m\times n}$, $A\in L^{1^*}(\Omega;\R^{m\times n})$ and a Borel set $\omega\subset\Omega$ we set
\[
\E_K(A;\omega):=\int_\omega \dist_K^{1^*}(A)\,\d x+|\curl A|(\omega)^{1^*} .
\]
We record two elementary properties of $\E_K$.
\begin{itemize}
	\item \emph{Scaling.} If $B(y):=A(x_0+ry)$ and $\omega_r:=r^{-1}(\omega-x_0)$, then, since $|\curl B|(\omega_r)=r^{1-n}|\curl A|(\omega)$ and $(1-n)1^*=-n$,
	\begin{equation}\label{e:scaling}
		\E_K(B;\omega_r)=r^{-n}\,\E_K(A;\omega).
	\end{equation}
	\item \emph{Superadditivity.} If $\{\omega_k\}$ is a countable family of pairwise disjoint Borel subsets of $\omega$, then, since $1^*\ge1$,
	\begin{gather*}
	\sum_k|\curl A|(\omega_k)^{1^*} \le\Big(\sum_k|\curl A|(\omega_k)\Big)^{1^*}\le|\curl A|(\omega)^{1^*},\\
	\text{so that}\qquad
	\sum_k\E_K(A;\omega_k)\le\E_K(A;\omega).
	\end{gather*}
\end{itemize}

The main estimate for incompatible fields is the following.

\begin{theorem}\label{t:rigidity}
	Let $\Omega \subset \R^n$ be a connected bounded open set with Lipschitz boundary, and let $K_1,K_2\subset \R^{m \times n}$ be two disjoint strongly incompatible compact sets, $K:=K_1\cup K_2$.
	Then there exists a constant $C=C(\Omega,K_1,K_2)>0$ such that, for every $A \in L^{1^*}(\Omega;\R^{m \times n})$,
	\begin{equation}\label{e:rigidity}
		\min \Big\{ \| \dist_{K_1} (A)\|_{L^{1^*}(\Omega)}, \| \dist_{K_2} (A)\|_{L^{1^*}(\Omega)} \Big\} \leq C \Big( \| \dist_{K} (A)\|_{L^{1^*}(\Omega)} + |\curl A|(\Omega) \Big).
	\end{equation}
\end{theorem}

By the elementary inequalities $(a+b)^{1^*}\le2^{1^*-1}(a^{1^*}+b^{1^*})$ and $a^{1^*}+b^{1^*}\le(a+b)^{1^*}$ for $a,b\ge0$, \eqref{e:rigidity} is equivalent, up to changing the value of $C$, to
\begin{equation}\label{e:rigidity2}
	\min_{i=1,2}\int_\Omega \dist_{K_i}^{1^*}(A)\,\d x\le C\,\E_K(A;\Omega).
\end{equation}

For the application to semi-coherent interfaces we only need $n=m=2$, so that $1^*=2$, and wells of the form $K_1=\SO(2)$ and $K_2=\mu \SO(2)$.
For multiple wells in two dimensions we recall the following result, due to \v{S}ver\'ak \cite{Sverak1993}:

\medskip

\emph{Let $U_1,\dots,U_N\in\R^{2\times2}$ with $\det U_i>0$, and assume that $K:=\bigcup_{i=1}^N\SO(2)U_i$ contains no rank-one connections, i.e.\ $\operatorname{rank}(X-Y)=2$ for all distinct $X,Y\in K$.
	Then, for every partition $\{1,\dots,N\}=I\cup J$ into disjoint nonempty sets, $\bigcup_{i\in I}\SO(2)U_i$ and $\bigcup_{i\in J}\SO(2)U_i$ are disjoint strongly incompatible sets.}

\medskip

In the case of interest for us, the strong incompatibility can be checked directly.

\begin{lemma}\label{l:wells}
	For every $\mu\in(0,1)$, the sets $\SO(2)$ and $\mu\SO(2)$ are disjoint strongly incompatible compact sets.
	More precisely, for every connected open set $\Omega\subset\R^2$ and every $\nu\in\GY^\infty(\Omega;\R^{2\times2})$ with $\supp\nu_x\subset K_\mu$ for a.e.\ $x\in\Omega$, there exists $S\in K_\mu$ such that $\nu_x=\delta_S$ for a.e.\ $x\in\Omega$.
\end{lemma}

\begin{proof}
	Let $\nu$ be generated by $(\nabla u_j)$, with $(u_j)$ bounded in $W^{1,\infty}(\Omega;\R^2)$ (up to subtracting constants). Up to a subsequence, $u_j\overset{*}{\rightharpoonup} u$ in $W^{1,\infty}$, and $\langle\nu_x,\mathrm{id}\rangle=\nabla u(x)$ for a.e.\ $x$. By the weak$^*$ continuity of the determinant on gradients bounded in $L^\infty$, we also have $\langle\nu_x,\det\rangle=\det\nabla u(x)$ for a.e.\ $x$.
	Every matrix $P\in K_\mu$ is conformal, so that $\det P=\frac12|P|^2$; moreover $\det M\le\frac12|M|^2$ for every $M\in\R^{2\times2}$. Hence, by Jensen's inequality, for a.e.\ $x$
	\[
	\frac12\langle\nu_x,|\cdot|^2\rangle=\langle\nu_x,\det\rangle=\det\nabla u(x)\le\frac12|\nabla u(x)|^2=\frac12\big|\langle\nu_x,\mathrm{id}\rangle\big|^2\le\frac12\langle\nu_x,|\cdot|^2\rangle .
	\]
	Equality holds in Jensen's inequality for the strictly convex function $|\cdot|^2$, so that $\nu_x=\delta_{\nabla u(x)}$, and $\nabla u(x)\in K_\mu$ for a.e.\ $x\in\Omega$.
	In particular $\nabla u$ is a conformal field with vanishing curl, i.e.\ a harmonic conformal field in the sense of Definition~\ref{d:harmonic conformal}, and Lemma~\ref{l:h_c-rigidity-1} gives that $\nabla u$ is constant in $\Omega$.
	Finally, $\SO(2)$ and $\mu\SO(2)$ are disjoint since $\mu\ne1$.
\end{proof}

Note that $K_\mu$ has no rank-one connections for $\mu\neq1$: for $R,R'\in\SO(2)$ one has $$\det(R-\mu R')=1+\mu^2-\mu\,R:R'\ge(1-\mu)^2>0.$$
A consequence of Theorem~\ref{t:rigidity} is the following rigidity estimate for incompatible fields and multiple wells.

\begin{corollary}\label{c:rigidity}
	Let $\Omega\subset\R^2$ be a simply connected bounded open set with Lipschitz boundary and let $K:=\bigcup_{i=1}^N\SO(2)U_i$ with $\det U_i>0$ and no rank-one connections. There exists $C=C(\Omega,K)$ such that for every $A\in L^2(\Omega;\R^{2\times2})$
	\begin{equation}\label{e:rigidity0}
		\inf_{S\in K}\|A-S\|_{L^2(\Omega)}\le C\Big(\|\dist_K(A)\|_{L^2(\Omega)}+|\curl A|(\Omega)\Big).
	\end{equation}
	In particular, by Lemma~\ref{l:wells}, \eqref{e:rigidity0} holds for $K=K_\mu$, $\mu\in(0,1)$, with a constant $C=C(\Omega,\mu)$.
\end{corollary}

\begin{proof}
	Here $n=2$ and $1^*=2$.

	\emph{Step 1: choice of the well.} We prove by induction on $N$ that
	\begin{equation}\label{e:choice well}
		\min_{i=1,\dots,N}\int_\Omega \dist^2_{\SO(2)U_i}(A)\,\d x\le C\,\E_K(A;\Omega).
	\end{equation}
	For $N=1$ there is nothing to prove. For $N\ge2$, let $K_1:=\SO(2)U_1$ and $K':=\bigcup_{i\ge2}\SO(2)U_i$, which are disjoint and strongly incompatible by the result of \v{S}ver\'ak recalled above. By \eqref{e:rigidity2}, either $\int_\Omega\dist^2_{K_1}(A)\le C\E_K(A;\Omega)$, and \eqref{e:choice well} follows, or $\int_\Omega\dist^2_{K'}(A)\le C\E_K(A;\Omega)$. In the second case, since $K'$ has no rank-one connections, the inductive hypothesis gives
	\[
	\min_{i\ge2}\int_\Omega \dist^2_{\SO(2)U_i}(A)\,\d x\le C\,\E_{K'}(A;\Omega)=C\Big(\int_\Omega\dist^2_{K'}(A)\,\d x+|\curl A|(\Omega)^2\Big)\le C\,\E_K(A;\Omega).
	\]

	\emph{Step 2: one-well rigidity.} Let $i$ realize the minimum in \eqref{e:choice well}, write $U:=U_i$, and set $\tilde A(y):=A(U^{-1}y)\,U^{-1}$ for $y\in U\Omega$, which is a connected bounded Lipschitz domain. The rows of $\tilde A$ are the pull-backs of the rows of $A$ under the linear map $y\mapsto U^{-1}y$, so that $\curl\tilde A=\det(U^{-1})\,(\curl A)\circ U^{-1}$ and, since $\det U>0$, $|\curl\tilde A|(U\Omega)=|\curl A|(\Omega)$. Moreover $$\dist_{\SO(2)}(\tilde A(y))\le|U^{-1}|\,\dist_{\SO(2)U}(A(U^{-1}y)).$$
	The one-well estimate of \cite{MullerScardiaZeppieri2014} on $U\Omega$ provides $R\in\SO(2)$ with
	\[
	\|\tilde A-R\|_{L^2(U\Omega)}\le C\Big(\|\dist_{\SO(2)}(\tilde A)\|_{L^2(U\Omega)}+|\curl\tilde A|(U\Omega)\Big).
	\]
	Changing variables back, $\|A-RU\|_{L^2(\Omega)}\le C\big(\|\dist_{\SO(2)U}(A)\|_{L^2(\Omega)}+|\curl A|(\Omega)\big)$, with $C$ depending also on $U$. Together with \eqref{e:choice well}, this gives \eqref{e:rigidity0} with $S:=RU\in K$.
\end{proof}

\subsection{Preparatory lemmas}
The following observation underlies the proof of Theorem~\ref{t:rigidity}.

\begin{proposition}\label{p:GY}
	Let $\Omega\subset \R^n$ be an open set, $K\subset\R^{m\times n}$ compact, and let $A_j\in L^{1^*}(\Omega;\R^{m\times n})$ satisfy
	\begin{equation}\label{e:vanishing}
		\lim_{j\to\infty}\left(\int_\Omega \dist_K^{1^*}(A_j)\,\d x+|\curl A_j|(\Omega)\right)=0 .
	\end{equation}
	Then, up to a subsequence, $(A_j)$ generates a Young measure $\nu$ such that
	\begin{enumerate}
		\item[(i)] $\supp\nu_x\subset K$ for a.e.\ $x\in\Omega$;
		\item[(ii)] for every open cube $Q\subset\Omega$, the restriction $\nu|_Q$ belongs to $\GY^\infty(Q;\R^{m\times n})$.
	\end{enumerate}
\end{proposition}

\begin{proof}
	By Lemma~\ref{l:pointwise}(a) and \eqref{e:vanishing}, $(A_j)$ is bounded in $L^{1^*}(\Omega)$, and by (Y3) a subsequence generates a Young measure $\nu$.
	By \eqref{e:vanishing}, $\dist_K^{1^*}(A_j)\to0$ in $L^1(\Omega)$; in particular this sequence is equi-integrable, and (Y1) gives
	\[
	\int_\Omega\langle\nu_x,\dist_K^{1^*}\rangle\,\d x=\lim_j\int_\Omega \dist_K^{1^*}(A_j)\,\d x=0,
	\]
	which is (i).

	For (ii), fix an open cube $Q\subset\Omega$. Let $s:=1^*$ if $n\ge3$ and $s\in(1,2)$ if $n=2$. By \cite{ContiGarroni2021} for $n\geq3$ and \cite{MullerScardiaZeppieri2014} for $n=2$, we can write $A_j=\nabla v_j+Y_j$ in $Q$, with $v_j\in W^{1,s}(Q;\R^m)$ and
	\[
	\|Y_j\|_{L^s(Q)}\leq C(n,s,Q)\,|\curl A_j|(Q) \to 0 .
	\]
	Hence $(\nabla v_j)$ is bounded in $L^s(Q)$ and, since $A_j-\nabla v_j\to0$ in measure, by (Y3) it generates the same Young measure $\nu|_Q$. Thus $\nu|_Q\in\GY^s(Q)$ and, being supported in the compact set $K$ by (i), $\nu|_Q\in\GY^\infty(Q)$ (see e.g.\ \cite{KinderlehrerPedregal1994}, \cite[Section~4]{Muller1999}).
\end{proof}

The following lemma shows that local gradient Young measures are enough to detect strongly incompatible wells.

\begin{lemma}\label{l:connected}
	Let $K_1,K_2$ be disjoint strongly incompatible compact sets, let $\Omega\subset\R^n$ be a connected open set and let $\nu$ be a Young measure on $\Omega$ with $\supp \nu_x\subset K_1\cup K_2$ for a.e.\ $x$ and $\nu|_Q\in\GY^\infty(Q)$ for every open cube $Q\subset\Omega$. Then there is $i\in\{1,2\}$ with $\supp\nu_x\subset K_i$ for a.e.\ $x\in\Omega$.
\end{lemma}

\begin{proof}
	By Definition~\ref{d:incompatible} applied in each cube $Q\subset\Omega$, there is $i(Q)\in\{1,2\}$ with $\supp\nu_x\subset K_{i(Q)}$ for a.e.\ $x\in Q$; $i(Q)$ is unique because $K_1\cap K_2=\emptyset$.
	Let $\Omega_i$ be the union of the open cubes $Q\subset\Omega$ with $i(Q)=i$. Then $\Omega_1,\Omega_2$ are open and cover $\Omega$. If $x\in\Omega_1\cap\Omega_2$, there are cubes $Q,Q'\ni x$ with $i(Q)=1$ and $i(Q')=2$; the open set $Q\cap Q'$ has positive measure, and on it $\nu_x$ would be a probability measure supported in $K_1\cap K_2=\emptyset$, which is impossible. Hence $\Omega_1\cap\Omega_2=\emptyset$ and, by connectedness, one of the two sets is empty.
\end{proof}

\subsection{Proof of Theorem~\ref{t:rigidity}}
We follow the blow-up argument of \cite{DeLellisSzekelyhidi2006} based on the transition layer estimate by Ball and James \cite{BallJames2015}, and we prove \eqref{e:rigidity2}.
Throughout the proof, $\eta$, $D$ and $C_K$ are the constants of Lemma~\ref{l:pointwise}, and $C$ denotes a constant depending only on $n$, $K_1$, $K_2$ and $\Omega$, which may change from line to line.
By Lemma~\ref{l:pointwise}(b),
\begin{equation}\label{e:trivial}
	\int_\omega\dist^{1^*}_{K_i}(A)\,\d x\le C\int_\omega\dist^{1^*}_{K}(A)\,\d x+C|\omega|\qquad\text{for every Borel }\omega\subset\Omega,\ i=1,2 .
\end{equation}

\medskip

\noindent\emph{Step 1: reduction to graph domains.}
Since $\Omega$ is a connected bounded Lipschitz domain, we can write $\Omega=\Omega_1\cup\dots\cup\Omega_\ell$, where each $\Omega_k$ is, up to a rigid motion, of the form
\begin{equation}\label{e:graph}
	\Omega_k=\Big\{(x',x_n):\ x'\in(0,r_0)^{n-1},\ f_k(x')<x_n<r_0\Big\},
\end{equation}
with $r_0>0$ and $f_k:[0,r_0]^{n-1}\to\R$ Lipschitz with $\max f_k<r_0$ (for interior pieces, $f_k$ is constant).
Since $\Omega$ is connected, the graph on $\{1,\dots,\ell\}$ in which $k$ and $k'$ are joined whenever $\Omega_k\cap\Omega_{k'}\neq\emptyset$ is connected. Set $\alpha:=\min\{|\Omega_k\cap\Omega_{k'}|:\ \Omega_k\cap\Omega_{k'}\neq\emptyset\}>0$.

Suppose that \eqref{e:rigidity2} holds in each $\Omega_k$ with constant $C_k$, and set $C_*:=\max_kC_k$. Given $A$, let $i_k$ be an index realizing the minimum in \eqref{e:rigidity2} on $\Omega_k$.
If $i_k=i$ for every $k$, then, by superadditivity,
\[
\int_\Omega \dist_{K_i}^{1^*}(A)\,\d x\le\sum_k C_*\,\E_K(A;\Omega_k)\le\ell\, C_*\,\E_K(A;\Omega).
\]
Otherwise, by the connectedness of the graph, there are $k,k'$ with $\Omega_k\cap\Omega_{k'}\ne\emptyset$ and $i_k\neq i_{k'}$. By Lemma~\ref{l:pointwise}(c),
\begin{align*}
2^{1-1^*}\eta^{1^*}\alpha&\le\int_{\Omega_k\cap\Omega_{k'}}\big(\dist_{K_1}^{1^*}(A)+\dist_{K_2}^{1^*}(A)\big)\, \d x\\
&\le\int_{\Omega_k}\dist^{1^*}_{K_{i_k}}(A)\,\d x+\int_{\Omega_{k'}}\dist^{1^*}_{K_{i_{k'}}}(A)\,\d x
\le 2C_*\,\E_K(A;\Omega),
\end{align*}
and therefore, by \eqref{e:trivial},
\[
\int_\Omega \dist^{1^*}_{K_1}(A)\,\d x\leq C \int_\Omega \dist^{1^*}_{K}(A)\, \d x + C|\Omega|
\leq \Big(C+\frac{2^{1^*}C\,C_*|\Omega|}{\eta^{1^*}\alpha} \Big)\, \E_K(A;\Omega).
\]
Hence it suffices to prove \eqref{e:rigidity2} for $\Omega$ of the form \eqref{e:graph}, which we assume from now on. We denote by $\Lambda$ the Lipschitz constant of the function $f$ defining $\Omega$ and set
\begin{equation}\label{e:L}
	L:=\max\Big\{1,\,4\sqrt{n-1}\,\Lambda\Big\},\qquad Q^L_1:=(-1,1)^{n-1}\times(-L,L),\qquad
	Q^L_r(x):=x+rQ^L_1 .
\end{equation}

\medskip

\noindent\emph{Step 2: contradiction assumption.}
If \eqref{e:rigidity2} fails, there are $A_j\in L^{1^*}(\Omega;\R^{m\times n})$ with $\E_K(A_j;\Omega)<\infty$ and
\begin{equation}\label{e:contra}
	\min_{i=1,2}\int_\Omega \dist_{K_i}^{1^*}(A_j)\,\d x>j\,\E_K(A_j;\Omega).
\end{equation}
By \eqref{e:trivial}, the left-hand side is at most $C_1\E_K(A_j;\Omega)+C_1|\Omega|$ for a suitable constant $C_1$, so that
\begin{equation}\label{e:Eto0}
	\E_K(A_j;\Omega)\le\frac{C_1|\Omega|}{j-C_1}\to0 .
\end{equation}
By Proposition~\ref{p:GY} and Lemma~\ref{l:connected}, up to a subsequence $(A_j)$ generates a Young measure $\nu$ supported, say, in $K_2$ for a.e.\ $x\in\Omega$. Since $\dist^{1^*}_{K_i}(A_j)\le C\dist^{1^*}_K(A_j)+C$ is equi-integrable, (Y1) and Lemma~\ref{l:pointwise}(c) give
\begin{gather}
	\int_\Omega \dist_{K_2}^{1^*}(A_j)\, \d x\to\int_\Omega\langle\nu_x,\dist_{K_2}^{1^*}\rangle\,\d x=0,\label{e:global_1}\\
	\int_\Omega \dist_{K_1}^{1^*}(A_j)\,\d x\to\int_\Omega\langle\nu_x,\dist_{K_1}^{1^*}\rangle\, \d x
	\geq \eta^{1^*}|\Omega|=:c_0 >0.\label{e:global_2}
\end{gather}
Let
\begin{align*}
S_j:=\Big\{x\in\Omega:\ \dist_{K_1}(A_j(x))\le \dist_{K_2}(A_j(x))\Big\},\\
f_j:=\chi_{S_j}\,\dist_{K_2}^{1^*}(A_j),\qquad
g_j:=\chi_{\Omega\setminus S_j}\,\dist_{K_1}^{1^*}(A_j).
\end{align*}
Thus $\dist_K=\dist_{K_1}$ on $S_j$ and $\dist_K=\dist_{K_2}$ on $\Omega\setminus S_j$; $f_j$ and $g_j$ measure the distance from the ``wrong'' well. By Lemma~\ref{l:pointwise}(c),
\begin{equation}\label{e:f lower}
	f_j\geq \tfrac12\big(\dist_{K_1}^{1^*}(A_j)+\dist_{K_2}^{1^*}(A_j)\big)\chi_{S_j}\geq 2^{-1^*}\eta^{1^*} \chi_{S_j}.
\end{equation}
Since $\int_\Omega\dist_{K_2}^{1^*}(A_j)=\int f_j+\int_{\Omega\setminus S_j}\dist_K^{1^*}(A_j)$, \eqref{e:contra} gives
\[
\int f_j\, \d x>(j-1)\,\E_K(A_j;\Omega)\ge0 ;
\]
in particular $|S_j|>0$. On the other hand, $\int f_j\le\int_\Omega \dist_{K_2}^{1^*}(A_j)\to0$ by \eqref{e:global_1}, whereas $\int(f_j+g_j)\ge\int_\Omega \dist_{K_1}^{1^*}(A_j)$ because $\dist_{K_1}\le\dist_{K_2}$ on $S_j$. By \eqref{e:global_2}, for $j$ large
\begin{equation}\label{e:negative}
	\int(f_j-g_j)\,\d x=2\int f_j\,\d x-\int(f_j+g_j)\,\d x\le-\frac{c_0}2 .
\end{equation}

\medskip

\noindent\emph{Step 3: choice of a good box.}
Fix $R_0>0$ such that $\Omega\subset Q^L_{R_0}(x)$ for every $x\in\Omega$. For $x\in\Omega$ the function
\[
H_x(r):=\frac{1}{|Q^L_r|}\int_{Q^L_r(x)}(f_j-g_j)\,\d y
\]
is continuous on $(0,\infty)$, and $H_x(R_0)<0$ by \eqref{e:negative}. For a.e.\ $x_0\in S_j$ (Lebesgue points of $f_j$ and $g_j$ with respect to the boxes $Q^L_r(x_0)$), $\lim_{r\to0}H_{x_0}(r)=f_j(x_0)>0$ by \eqref{e:f lower}. Hence there is $r(x_0)\in(0,R_0)$ with
\[
\int_{Q^L_{r(x_0)}(x_0)}f_j\,\d x=\int_{Q^L_{r(x_0)}(x_0)}g_j\,\d x .
\]
By the Besicovitch covering theorem (in the version for homothetic copies of a fixed box centered at the points), there are $N_B=N_B(n,L)$ countable families $\mathcal Q_1,\dots,\mathcal Q_{N_B}$ of pairwise disjoint such boxes whose union covers almost every point of $S_j$, hence the support of $f_j$ up to a null set.
For some $\bar k$, by superadditivity,
\[
\sum_{Q\in\mathcal Q_{\bar k}}\int_Qf_j\,\d x\ge\frac1{N_B}\int f_j\,\d x
>\frac{j-1}{N_B}\,\E_K(A_j;\Omega)
\ge\frac{j-1}{N_B}\sum_{Q\in\mathcal Q_{\bar k}}\E_K(A_j;Q\cap\Omega).
\]
Therefore there is a box $\bar Q=Q^L_{r_j}(x_j)\in\mathcal Q_{\bar k}$, with $x_j\in S_j$ and $r_j<R_0$, such that
\begin{equation}\label{e:goodbox}
	\int_{\bar Q}g_j\,\d x=\int_{\bar Q}f_j\,\d x
	>\frac{j-1}{N_B}\,\E_K(A_j;\bar Q\cap\Omega).
\end{equation}

\medskip

\noindent\emph{Step 4: blow-up.}
Set
\[
U_j:=r_j^{-1}(\Omega-x_j)\cap Q^L_1,\qquad
\Sigma_j:=r_j^{-1}(S_j-x_j)\cap Q^L_1\subset U_j,
\qquad B_j(y):=A_j(x_j+r_jy),\ y\in U_j .
\]
Since $f_j$, $g_j$ and $\E_K$ all scale like $r_j^{-n}$ by \eqref{e:scaling}, \eqref{e:goodbox} gives
\begin{equation}\label{e:rescaled}
	\int_{\Sigma_j}\dist_{K_2}^{1^*}(B_j)\, \d y=\int_{U_j\setminus\Sigma_j}\dist_{K_1}^{1^*}(B_j)\,\d y
	>\frac{j-1}{N_B}\,\E_K(B_j;U_j).
\end{equation}
By \eqref{e:trivial}, the left-hand side of \eqref{e:rescaled} is at most $C_1\E_K(B_j;U_j)+C_1|Q^L_1|$, and therefore, as in \eqref{e:Eto0}, $\E_K(B_j;U_j)\to0$.

By Lemma~\ref{l:C}, applied with $x_j\in\Omega$ and $r_j\in(0,R_0]$, we have $|U_j|\ge\gamma>0$ and, up to a subsequence, there is a connected open set $U\subset Q^L_1$ with $|U_j\triangle U|\to0$ (so that $|U|\ge\gamma$). Moreover, for every $\delta \in (0, \gamma)$, to be chosen later, there is a connected open set $U'\subset U$ with $|U\setminus U'|<\delta$ and $U'\subset U_j$ for all $j$ large.

\medskip

\noindent\emph{Step 5: conclusion.}
Since $\E_K(B_j;U')\le\E_K(B_j;U_j)\to0$, Proposition~\ref{p:GY} and Lemma~\ref{l:connected} applied in the connected set $U'$ give, up to a subsequence, either $\int_{U'}\dist^{1^*}_{K_1}(B_j)\,\d y\to0$ or $\int_{U'}\dist^{1^*}_{K_2}(B_j)\,\d y\to0$. We show that both alternatives lead to a contradiction. Set
\[
a_j:=\int_{\Sigma_j\cap U'}\dist^{1^*}_{K_2}(B_j)\, \d y-\int_{U'\setminus\Sigma_j}\dist^{1^*}_{K_1}(B_j)\, \d y.
\]
By the equality in \eqref{e:rescaled}, by \eqref{e:trivial} and since $U'\subset U_j$ for $j$ large,
\begin{align}\label{e:a_j}
	|a_j| & = \left\vert
	\int_{\Sigma_j\cap (U_j\setminus U')}\dist^{1^*}_{K_2}(B_j) \,\d y-\int_{(U_j\setminus U')\setminus\Sigma_j}\dist^{1^*}_{K_1}(B_j) \,\d y
	\right\vert \notag\\
	& \leq C \int_{U_j} \dist_K^{1^*}(B_j) \,\d y + C\,|U_j\setminus U'|\leq C\, \E_K(B_j; U_j) + C\,|U_j\triangle U| + C \delta,
\end{align}
so that $\limsup_j|a_j|\le C\delta$.

If $\int_{U'}\dist^{1^*}_{K_1}(B_j) \,\d y\to0$, using that $\dist_{K_2}(B_j)\leq\dist_{K_1}(B_j)$ in $U'\setminus \Sigma_j$ we get
\begin{align*}
	\limsup_j \int_{U'}\dist^{1^*}_{K_2}(B_j) \,\d y & = \limsup_j\Big( a_j +  \int_{U'\setminus \Sigma_j}\dist^{1^*}_{K_1}(B_j) \,\d y + \int_{U'\setminus \Sigma_j}\dist^{1^*}_{K_2}(B_j) \,\d y\Big)\\
	&\le C\delta +  2 \limsup_j \int_{U'\setminus \Sigma_j}\dist^{1^*}_{K_1}(B_j) \,\d y = C\delta.
\end{align*}
Similarly, if $\int_{U'}\dist^{1^*}_{K_2}(B_j) \,\d y\to0$, using that $\dist_{K_1}(B_j)\leq\dist_{K_2}(B_j)$ in $U'\cap\Sigma_j$ we get
\begin{align*}
	\limsup_j \int_{U'}\dist^{1^*}_{K_1}(B_j) \,\d y
	& = \limsup_j\Big( - a_j +  \int_{U'\cap \Sigma_j}\dist^{1^*}_{K_2}(B_j) \,\d y + \int_{U'\cap \Sigma_j}\dist^{1^*}_{K_1}(B_j) \,\d y\Big)\\
	&\le C\delta +  2 \limsup_j \int_{U'\cap \Sigma_j}\dist^{1^*}_{K_2}(B_j) \,\d y = C\delta.
\end{align*}
In both cases, by Lemma~\ref{l:pointwise}(c) and $|U'|\ge|U|-\delta\ge\gamma-\delta$,
\[
2^{1-1^*}\eta^{1^*}(\gamma-\delta)\le2^{1-1^*}\eta^{1^*}|U'|
\le\limsup_{j}\int_{U'}\big(\dist^{1^*}_{K_1}(B_j)+\dist^{1^*}_{K_2}(B_j) \big) \,\d y
\le C\delta,
\]
which is a contradiction if $\delta$ is chosen sufficiently small. This concludes the proof of Theorem~\ref{t:rigidity}. \qed

\medskip

We conclude this section by proving the compactness result used in the proof
of Theorem~\ref{t:main}. The argument follows that of
\cite[Lemma 4.6]{FortunaGarroniSpadaro2025}.

\begin{proposition}\label{p:compactness}
	Let $\Omega \subset \R^2$ be a bounded open set, $K:=\bigcup_{i=1}^N\SO(2)U_i$ with $\det U_i>0$ and no rank-one
	connections, and $A_j \in L^2(\Omega;\R^{2 \times 2})$ such that
	\begin{equation}\label{e:comp}
		\lim_{j \to \infty} \|\dist_K(A_j)\|_{L^2(\Omega)}=0, \qquad \sup_{j \in \N} |\curl(A_j)|(\Omega)< \infty.
	\end{equation}
	Then, up to a subsequence, $A_j$ converges to some $A\in BV(\Omega;K)$ weakly in $L^2$ and strongly in $L^1$. Moreover,
	\begin{equation}\label{e:curlcontrolsgradient}
		|DA| \leq C|\curl A|.
	\end{equation} 
\end{proposition}

\begin{proof}
	We begin by recalling the following fact, see Corollary 4.3
	in \cite{FortunaGarroniSpadaro2025}.
	Let $(A_j)_j \subset L^1(\Omega;\R^{2\times 2})$ be a sequence
	converging weakly to $A\in L^1(\Omega;\R^{2\times 2})$,
	let $(\lambda_j)_j$ be a sequence of positive measures
	converging weakly$^*$ to $\lambda$, and let $(\gamma_j)_j$
	be a sequence converging to zero.
	Assume that, for every cube $Q\subset \Omega$, there exists
	a constant matrix $A_{j,Q}\in K$ such that
	\[
	\int_Q |A_j(x)-A_{j,Q}|\,\d x
	\leq \ell(Q)\lambda_j(Q)+\gamma_j,
	\]
	where $\ell(Q)$ denotes the side length of $Q$.
	Then $A_j\to A$ strongly in $L^1(\Omega;\R^{2\times 2})$,
	$A\in BV(\Omega;K)$, and there exists a constant $C>0$
	such that $|DA|\leq C\lambda$.
	
	Now consider a sequence $(A_j)_j$ satisfying \eqref{e:comp}.
	It follows from Lemma \ref{l:pointwise} that
	$(A_j)_j$ is bounded in $L^2(\Omega;\R^{2\times 2})$.
	Hence there exists $A\in L^2(\Omega;\R^{2\times 2})$ such that,
	up to a subsequence, not relabeled, $A_j$ converges weakly
	to $A$ in $L^2(\Omega;\R^{2\times 2})$.
	Moreover, by the uniform boundedness of the curls,
	we may extract a further subsequence such that
	$\curl A_j$ converges weakly$^*$ to $\curl A$.
	
	For every cube $Q\subset \Omega$, we can apply Corollary \ref{c:rigidity} to find $A_{j,Q}\in K$ such that
	\[
	\|A_j-A_{j,Q}\|_{L^2(Q)}
	\leq C\left(
	\|\dist_K(A_j)\|_{L^2(Q)}
	+|\curl(A_j)|(Q)
	\right).
	\]
	By H\"older's inequality,
	\begin{align*}
		\|A_j-A_{j,Q}\|_{L^1(Q)}
		&\leq \ell(Q)\|A_j-A_{j,Q}\|_{L^2(Q)}\\
		&\leq C\ell(Q)\left(
		\|\dist_K(A_j)\|_{L^2(Q)}
		+|\curl(A_j)|(Q)
		\right)\\
		&\leq C\ell(Q)|\curl(A_j)|(Q)
		+C\ell(Q_L)\|\dist_K(A_j)\|_{L^2(Q_L)}.
	\end{align*}
	Choosing
	\[
	\gamma_j:=L\|\dist_K(A_j)\|_{L^2(\Omega)},
	\qquad
	\lambda_j:=C|\curl(A_j)|,
	\]
	we conclude by the compactness criterion recalled above.
	To obtain \eqref{e:curlcontrolsgradient}, we apply the rigidity
	estimate and H\"older's inequality to $A$ on each cube
	$Q\subset \Omega$. Since $A$ takes values in $K$ (i.e. $\dist_K(A)=0$), we obtain
	\[
	\|A-A_Q\|_{L^1(Q)} \leq C\ell(Q)|\curl(A)|(Q)
	\]
	for some constant matrix $A_Q\in K$.
	Thus, $A$ satisfies a Poincar\'e inequality with respect to
	$|\curl A|$, which controls $|DA|$ by a standard
	argument.
\end{proof}

\section{Energy transfer estimate}\label{s:transfer}
\subsection{Regularization and harmonic replacement}

In this subsection we recall the approximation results by smoother competitors
for sequences $(A_\eps,S_\eps)\in\mathcal{C}_\eps(\mu)$ in the interfacial energy scaling regime.
The arguments closely follow those developed by Lauteri--Luckhaus
\cite{LauteriLuckhaus2016} and revisited in
\cite{FortunaGarroniSpadaro2025}, with minor adaptations.
For the reader's convenience, we state the results and briefly explain
the differences from these references.
These approximation results do not rely on the prescribed boundary data; in fact, the modified fields constructed below do not satisfy \eqref{e:prescribed-phase} in general. For this reason we introduce the larger class
\[
\Chat_\eps\supset\mathcal C_\eps(\mu)
\]
of all configurations $(A,S)$, with $S\subset\R^2$ closed, $A\in L^1_{\rm loc}(\R^2;\R^{2\times2})$ and $\spt(\curl A)\subset S$, which satisfy the quantization condition \eqref{e:quantization}, with no condition on the boundary values. The boundary condition \eqref{e:prescribed-phase} enters the proof of Theorem~\ref{t:main} only through the original configurations $(A_j,S_j)\in\mathcal C_{\eps_j}(\mu)$, see \eqref{e:phases}.

The first lemma establishes an energetic equivalence between general competitors and regular competitors with controlled curl.

\begin{lemma}\label{l:asymptoticequivalent}
	There exists a universal constant $C>0$ with the following property. For every $L>3\eps$ and every $(A,S)\in\Chat_{\varepsilon}$, there exists
	$(\widetilde{A},\widetilde{S})\in \Chat_{\varepsilon}$ such that for every $0<L'\leq L-3\varepsilon$
	\begin{itemize}
		\item[(i)]	\(
		\|\widetilde{A}\|_{L^{\infty}(Q_{L'})}\leq C\);
		\item[(ii)]
		\(
		\left\|
		\widetilde{A}
		-A\,\mathbbm{1}_{Q_{L'}\setminus B_{\lambda\varepsilon}(S)}
		\right\|^2_{L^2(Q_{L'})}
		+
		\mathcal{F}_{\varepsilon}(\widetilde{A},\widetilde{S};Q_{L'})
		\leq
		C \mathcal{F}_{\varepsilon}(A,S;Q_L);
		\)
		\item[(iii)]
		$\curl\widetilde{A}\in L^{\infty}(Q_{L'};\R^2)$ and
		\(
		|\curl\widetilde{A}(x)|
		\leq
		\frac{C}{\lambda\varepsilon}
		\mathbbm{1}_{\widetilde{S}}(x)\)
		for every  $x\in Q_{L'}$.
	\end{itemize}
\end{lemma}

We omit the proof, which essentially follows from combining
the proofs of Lemmas 3.1 and 3.2 in \cite{FortunaGarroniSpadaro2025}, together with the observation that the specific structure of the wells is never used in the arguments:
the same proof remains valid if 
$\operatorname{SO}(2)$ is replaced by any compact set $K\subset\R^{2\times 2}$ in \cite{FortunaGarroniSpadaro2025}.
We also observe that the construction preserves the quantization condition \eqref{e:quantization}, but it is not designed to preserve the exact values $\Id$ and $\mu\Id$ of $A$ on the half-planes $\{\pm x_2>\delta\}$ (the modifications are local, at scale $\eps$, and may reach these half-planes from the strip). Hence \eqref{e:prescribed-phase} is in general not preserved, and the statement is formulated in the class $\Chat_\eps$.

\medskip

The next lemma establishes the existence of a harmonic competitor with comparable energy.
The proof closely follows that of \cite{LauteriLuckhaus2016}.
The main difference is the presence of two wells, whose conformal structure
is essential to the argument. We use the cofactor estimate
\eqref{e:projection} and the harmonic replacement argument of
Lemma~\ref{l:harmonic_conformal} on the defect-free part of the smaller cube.
For more general wells satisfying the assumptions of
Corollary~\ref{c:rigidity}, the corresponding equation may be quasilinear
rather than linear, so the same argument does not apply directly;
see, for instance, \cite{ChaudhuriMuller2003}.


\begin{lemma}\label{l:harmonic}
	There exists a universal constant $C>0$ with the following property. For every
	$(A,S)\in\Chat_{\varepsilon}$ and every pair of cubes
	$Q_L,Q_{L'}$ with $0<L'\leq L-3\varepsilon$, there exists
	$(A_{\mathrm{harm}},S_{\mathrm{harm}})
	\in\Chat_{\varepsilon}$ such that
	\begin{itemize}
		\item[(i)]
		$A_{\mathrm{harm}}$ is harmonic in
		$Q_{L'}\setminus S_{\mathrm{harm}}$ ;
		
		\item[(ii)]	\( 
		\left\|
		A_{\mathrm{harm}}
		-A\,\mathbbm{1}_{Q_{L'}\setminus B_{\lambda\varepsilon}(S)}
		\right\|^2_{L^2(Q_{L'})}
		+
		\mathcal{F}_{\varepsilon}(A_{\mathrm{harm}},S_{\mathrm{harm}};Q_{L'})
		\leq
		C \mathcal{F}_{\varepsilon}(A,S;Q_L);
		\)
		
		\item[(iii)]
		$\curl A_{\mathrm{harm}}\in L^{\infty}(Q_{L'};\R^2)$ and
		\(
		|\curl A_{\mathrm{harm}}(x)|
		\leq
		\frac{C}{\lambda\varepsilon}
		\mathbbm{1}_{S_{\mathrm{harm}}}(x)\)
		for every $x\in Q_{L'}$.		
	\end{itemize}
\end{lemma}
\begin{proof}
Let $(\widetilde A,\widetilde S)$ be the pair given by
Lemma~\ref{l:asymptoticequivalent}. 
Arguing as in Lemma~\ref{l:harmonic_conformal}, 
there exists a unique
$u\in H^1_0(Q_{L'}\setminus\widetilde S;\R^2)$ solving
\begin{equation}\label{e:pde}
    \Delta u=\operatorname{div}\bigl(\operatorname{cof}\widetilde A-\widetilde A\bigr)
    \qquad\text{in }Q_{L'}\setminus\widetilde S.
\end{equation}
Extend $u$ by zero to $\R^2$ and define
\begin{equation}\label{e:defAharm}
    A_{\mathrm{harm}}:=\widetilde A+\nabla u,
    \qquad S_{\mathrm{harm}}:=\widetilde S.
\end{equation}
Since $\curl\widetilde A=0$ in $Q_{L'}\setminus\widetilde S$, we have
$\operatorname{div}(\operatorname{cof}\widetilde A)=0$ there. Thus
\[
    \operatorname{div}A_{\mathrm{harm}}
    =\operatorname{div}\widetilde A+\Delta u=0,
    \qquad
    \curl A_{\mathrm{harm}}=\curl\widetilde A=0
    \quad\text{in }Q_{L'}\setminus\widetilde S,
\]
and consequently $\Delta A_{\mathrm{harm}}=0$ in this set, proving (i).
Moreover,
   $ \curl A_{\mathrm{harm}}=\curl\widetilde A$ in $\R^2$, because the extension of $u$ by zero belongs to $H^1(\R^2;\R^2)$ and therefore $\curl\nabla u=0$ in the sense of distributions in $\R^2$. In particular $\spt(\curl A_{\mathrm{harm}})\subset S_{\mathrm{harm}}$ and, since the circulation of $\nabla u$ along every closed rectifiable curve in $\R^2\setminus\widetilde S$ vanishes (in the same weak sense as in \eqref{e:quantization}), $A_{\mathrm{harm}}$ and $\widetilde A$ have the same circulations. Hence
$(A_{\mathrm{harm}},S_{\mathrm{harm}})\in\Chat_\varepsilon$,
and (iii) follows from Lemma~\ref{l:asymptoticequivalent}(iii).
Note that, even if $(A,S)\in\mathcal C_\eps(\mu)$, in general $(A_{\mathrm{harm}},S_{\mathrm{harm}})\notin\mathcal C_\eps(\mu)$: the function $u$ solves a Dirichlet problem in the whole of $Q_{L'}\setminus\widetilde S$ and does not vanish, in general, in $Q_{L'}\cap\{|x_2|>\delta\}$.

Testing \eqref{e:pde} with $u$ and using \eqref{e:projection} gives
\begin{equation*}
    \|\nabla u\|_{L^2(Q_{L'}\setminus\widetilde S)}
\leq 2\|\dist_{K_\mu}(\widetilde A)\|_{L^2(Q_{L'}\setminus\widetilde S)}.
\end{equation*}
Since $\widetilde A$ is bounded by Lemma~\ref{l:asymptoticequivalent}(i),
\begin{align*}
    \|A_{\mathrm{harm}}-\widetilde A\|_{L^2(Q_{L'})}^2
    &\le 4\int_{Q_{L'}\setminus\widetilde S}
       \dist_{K_\mu}^2(\widetilde A)\,\d x\\
    &\le 4\int_{Q_{L'}\setminus B_{\lambda\varepsilon}(\widetilde S)}
       \dist_{K_\mu}^2(\widetilde A)\,\d x
       +C|B_{\lambda\varepsilon}(\widetilde S)\cap Q_{L'}|\\
    &\le C\mathcal{F}_{\varepsilon}(\widetilde A,\widetilde S;Q_{L'}).
\end{align*}
Finally,
\[
    \mathcal{F}_{\varepsilon}(A_{\mathrm{harm}},\widetilde S;Q_{L'})
    \le 2\|A_{\mathrm{harm}}-\widetilde A\|_{L^2(Q_{L'})}^2
       +2\mathcal{F}_{\varepsilon}(\widetilde A,\widetilde S;Q_{L'}),
\]
and
\begin{align*}
    \|A_{\mathrm{harm}}-A\mathbbm{1}_{Q_{L'}\setminus B_{\lambda\varepsilon}(S)}\|_{L^2(Q_{L'})}^2
    &\le 2\|A_{\mathrm{harm}}-\widetilde A\|_{L^2(Q_{L'})}^2+2\|\widetilde A-A\mathbbm{1}_{Q_{L'}\setminus B_{\lambda\varepsilon}(S)}\|_{L^2(Q_{L'})}^2.
\end{align*}
Combining these inequalities with Lemma~\ref{l:asymptoticequivalent}(ii)
proves (ii).
\end{proof}

\subsection{Lauteri--Luckhaus foliation}
\begin{lemma}[Foliation]\label{l:foliation}
	There exist absolute constants $\delta_0\in(0,1)$ and $C>0$ with the following property.
	Let $R>0$, $x_0\in \R^2$ and $\{B_{\rho_i}(x_i)\}_{i=1,\dots,N}$ be a finite family of balls in $\mathbb R^2$ intersecting $B_{2R}\setminus \overline{B_R}(x_0)$ such that
	\begin{equation}\label{e:small radii}
		\sum_{i=1}^N\rho_i<\delta_0 \, R.
	\end{equation}
	Then there exists a Lipschitz function $\varphi\colon B_{2R}(x_0)\to[0,1]$ such that
	\begin{enumerate}
		\item $\|\nabla\varphi\|_\infty\le \frac{C}{R}$;
		\item $\varphi=1$ on $B_{\frac{9R}{8}}(x_0)$ and $\varphi=0$ in a neighborhood of $B_{2R}\setminus \overline{B_{\frac{15R}{8}}}(x_0)$;
		\item $\varphi$ is constant on a neighborhood of each ball $B_{\rho_i}(x_i)$;
		\item setting $U:=\big(B_{2R}\setminus \overline{B_R}(x_0)\big)\setminus\bigcup_{i=1}^N\overline{B_{\rho_i}(x_i)}$ and
		\[
		N':=\#\mathcal I',\qquad \mathcal I':=\Big\{i:\ B_{\rho_i}(x_i)\subset B_{2R}\setminus\overline{B_R}(x_0)\Big\},
		\]
		we have
		\begin{equation}\label{e:foliation estimate}
			\int_U\frac{|x-x_0|^2\,|\nabla\varphi(x)|^2}{\dist^2(x,\partial U)}\,dx\le C (1+N').
		\end{equation}
	\end{enumerate}
\end{lemma}

\begin{proof}
	We fix $M:=40$, $\kappa:=\frac M4-1=9$, $c_*:=\frac1{100}$ and $\delta_0:=\frac{c_*}{M}=\frac1{4000}$.	
	
	The lemma is invariant under translation and dilation; therefore we can assume that $x_0=0$ and $R=1$.
	Let $i\notin\mathcal I'$.
	Then $B_{\rho_i}(x_i)$ meets $\partial B_1$, or meets $\partial B_2$.
	Since $2\rho_i<2\delta_0$ by \eqref{e:small radii}, in all cases
	\begin{equation}\label{e:boundary balls}
		B_{\rho_i}(x_i)\subseteq B_{1+2\delta_0}\cup\big(\R^2\setminus \overline{B_{2-2\delta_0}}\big)\qquad\forall\; i\notin\mathcal I'.
	\end{equation}
	First we assume $N'\ge1$, and we set $I:=\{x_i:\ i\in\mathcal I'\}$, $d(x):=\dist(x,I)$,
	\[
	r_0:=\frac{\delta_0}{N'},\qquad r_k:=M^k r_0\quad (k\in\mathbb N),
	\]
	and we let $K\ge1$ be the largest integer such that $r_K\le c_*$ (note that $r_1\le M\delta_0=c_*$).
	In particular $r_{K-1}>c_*M^{-2}$.
	
	We now construct a family of nested nets in $I$.
		Let $I_{-1}:=I$ and, for $k=0,\dots,K$, let $I_k\subseteq I_{k-1}$ be a maximal subset with $|y-y'|\ge r_k$ for all distinct $y,y'\in I_k$.
	By maximality, every point of $I_{k-1}$ is at distance $<r_k$ from $I_k$.
	For $k=-1,\dots,K-1$ let $E_k\colon I_k\to I_{k+1}$ map each $w$ to a closest point of $I_{k+1}$, with $E_k(w)=w$ if $w\in I_{k+1}$; thus $|E_k(w)-w|<r_{k+1}$.
	For $z\in I_{k+1}$, its \emph{children} are the points of $\ch(z):=E_k^{-1}(z)\subseteq I_k$; note that $z\in\ch(z)$.
	For $x_i\in I$ we define its ancestors $y_{-1}(x_i):=x_i$ and $y_h(x_i):=E_{h-1}(y_{h-1}(x_i))\in I_h$.
	Then
	\begin{equation}\label{e:ancestors}
		|x_i-y_h(x_i)|<\sum_{m=0}^h r_m\le\frac{M}{M-1}\,r_h<2r_h .
	\end{equation}
	Since the sets $\ch(z)$, $z\in I_{k+1}$, partition $I_k$, the set
	\[
	D_{k+1}:=\{z\in I_{k+1}:\ \#\ch(z)\ge2\}
	\quad\text{satisfies}\quad
	|I_k|=\sum_{z\in I_{k+1}}\#\ch(z)\ge|I_{k+1}|+\#D_{k+1}.
	\]
	Summing over $k$, the sum telescopes:
	\begin{equation}\label{e:telescopic}
		\sum_{k=0}^{K-1}\#D_{k+1}\le|I_0|-|I_K|\le N' .
	\end{equation}
	
	Next we show that vertices with only one child give empty annuli.
		We say that a pair $(y,s)$ with $y\in I$ and $s>0$ is \emph{isolated} if
	\[
	I\cap\{x:\ s\le|x-y|\le\kappa s\}=\emptyset .
	\]
	We claim that if $z\in I_{k+1}\setminus D_{k+1}$, i.e.\ $\ch(z)=\{z\}$, then $(z,2r_k)$ is isolated.
	
	Let $w\in I_k$ with $|w-z|<r_{k+1}/2$.
	For every $z'\in I_{k+1}\setminus\{z\}$ we have
	$|w-z'|\ge|z-z'|-|w-z|>r_{k+1}/2>|w-z|$, hence $E_k(w)=z$, i.e.\ $w\in\ch(z)$, and therefore $w=z$.
	So $|w-z|\ge r_{k+1}/2$ for all $w\in I_k\setminus\{z\}$.
	Now let $x_m\in I$ and $w:=y_k(x_m)$.
	If $w=z$, then $|x_m-z|<2r_k$ by \eqref{e:ancestors}.
	Otherwise
	$|x_m-z|>\frac{r_{k+1}}2-2r_k=\big(\tfrac M2-2\big)r_k=2\kappa r_k$.
	This proves the claim.
	
	Isolated pairs are either nested or far apart.
		Let $(y,s)$ and $(y',t)$ be isolated, with $t\ge s$. We claim that
	\begin{equation}\label{e:dichotomy}
		\text{either}\quad |y-y'|<t\ \text{ and }\ B_{2s}(y)\subseteq B_{3t}(y'),
		\qquad\text{or}\quad |y-y'|>\kappa t .
	\end{equation}
	Since $y\in I$ and $(y',t)$ is isolated, either $|y-y'|<t$ or $|y-y'|>\kappa t$.
	In the first case, since $y'\in I$ and $(y,s)$ is isolated, either $|y-y'|<s$, and then $B_{2s}(y)\subset B_{3s}(y')\subset B_{3t}(y')$;
	or $\kappa s<|y-y'|<t$, and then $s<t/\kappa$ and $B_{2s}(y)\subset B_{t+2t/\kappa}(y')\subset B_{3t}(y')$.
	
	Consider the finite family
	\[
	\mathcal A:=\big\{(z,2r_k):\ 1\le k\le K-1,\ z\in I_{k+1}\setminus D_{k+1}\big\},
	\]
	whose elements are isolated by the previous claim.
	We select a subfamily $\mathcal M\subseteq\mathcal A$ greedily.
	Order $\mathcal A$ by nonincreasing radius.
	Select the first element $(y',t)$, discard all $(y,s)\in\mathcal A$ with $|y-y'|<t$, and repeat with the remaining elements.
	By \eqref{e:dichotomy}, the selected pairs $\mathcal M=\{(y_j,t_j)\}_j$ satisfy:
	\begin{enumerate}
		\item[(i)] $|y_j-y_l|>\kappa\max\{t_j,t_l\}$ for $j\ne l$. In particular the centers are distinct, so $\#\mathcal M\le N'$, and the balls $B(y_j,4t_j)$ are pairwise disjoint because $8<\kappa$.
		\item[(ii)] $B_{2s}(y)\subseteq G:=\bigcup_j B_{3t_j}(y_j)$ for every $(y,s)\in\mathcal A$.
	\end{enumerate}
	Moreover $4t_j\le 8r_{K-1}\le 8c_*/M<\frac1{16}$.
	
	We can now construct $\varphi$.
		Let $\varphi_0(x):=\min\{1,\max\{0,2(\tfrac74-|x|)\}\}$, so that $\varphi_0=1$ on $B_{\tfrac54}$, $\varphi_0=0$ outside $B_{\tfrac74}$ and $\Lip(\varphi_0)=2$.
	Let $\eta\colon[0,\infty)\to[0,1]$ be equal to $1$ on $[0,3]$, to $4-\tau$ on $[3,4]$ and to $0$ on $[4,\infty)$.
	Set $\eta_j(x):=\eta(|x-y_j|/t_j)$, let $m_j$ be the average of $\varphi_0$ on $B_{4t_j}(y_j)$, and define
	\[
	\varphi_1:=\sum_j\eta_j\,m_j+\Big(1-\sum_j\eta_j\Big)\varphi_0 .
	\]
	Since the supports of the $\eta_j$ are disjoint, on $B_{4t_j}(y_j)$ we have
	$\nabla\varphi_1=\nabla\eta_j\,(m_j-\varphi_0)+(1-\eta_j)\nabla\varphi_0$
	and $|m_j-\varphi_0|\le 8t_j\Lip(\varphi_0)$ there.
	Hence $\varphi_1$ takes values in $[0,1]$, $\Lip(\varphi_1)\le 9\Lip(\varphi_0)=18$, and $\varphi_1$ is constant on each $B_{3t_j}(y_j)$, i.e.\ on $G$.
	Since $4t_j<\frac1{16}$, every ball $B_{4t_j}(y_j)$ meeting $B_{1+\frac18}$ is contained in $B_{\frac54}$, and every such ball meeting $\{|x|\ge 2-\frac18\}$ is contained in $\{|x|\ge\frac74\}$.
	Therefore $\varphi_1=1$ on $B(0,1+\frac18)$ and $\varphi_1=0$ on $B_2\setminus B_{2-\frac18}$.
	
	Next, for $i\in\mathcal I'$, let $V_i:=B_{2\rho_i+r_0}(x_i)$, and let $\mathcal J:=\varphi_1\big(\bigcup_{i\in\mathcal I'}V_i\big)\subseteq[0,1]$.
	Each $\varphi_1(V_i)$ is an interval of length at most $2\Lip(\varphi_1)(2\rho_i+r_0)$, so by \eqref{e:small radii}
	\[
	|\mathcal J|\le 36\Big(2\sum_{i\in\mathcal I'}\rho_i+N'r_0\Big)\le 108\,\delta_0<\tfrac12 .
	\]
	Let $\psi(\tau):=\frac1{1-|\mathcal J|}\int_0^\tau\chi_{[0,1]\setminus\mathcal J}$.
	Then $\psi(0)=0$, $\psi(1)=1$, $\Lip(\psi)\le2$, and $\psi$ is constant on each interval $\varphi_1(V_i)$.
	We set $\varphi:=\psi\circ\varphi_1$.
	Then $\Lip(\varphi)\le36$, which gives (1).
	Moreover $\varphi=1$ on $B_{1+\frac18}$ and $\varphi=0$ on $B_2\setminus B_{2-\frac18}$, which gives (2).
	Finally $\varphi$ is constant on each ball $B_{3t_j}(y_j)$ and on each $V_i\supseteq B_{\rho_i}(x_i)$ with $i\in\mathcal I'$.
	For $i\notin\mathcal I'$, \eqref{e:boundary balls} gives $B(x_i,\rho_i)\subseteq B_{1+\frac18}\cup\{|x|>2-\frac18\}$, since $2\delta_0<\frac18$.
	Each such ball is connected, so it lies in one of the two sets, on which $\varphi$ equals $1$ and $0$ respectively.
	This gives (3).
	
	It remains to prove \eqref{e:foliation estimate}.
		Since $\nabla\varphi=0$ a.e.\ on $G\cup\bigcup_{i\in\mathcal I'}V_i\cup B_{1+\frac18}\cup\{|x|\ge2-\frac18\}$, it suffices to integrate over the set
	\[
	F:=\big\{x\in U:\ 1+\tfrac18<|x|<2-\tfrac18\big\}\setminus\Big(G\cup\bigcup_{i\in\mathcal I'}V_i\Big).
	\]
	For $x\in F$ and $i\in\mathcal I'$ we have $|x-x_i|\ge2\rho_i+r_0$, hence $\dist(x,\overline{B_{\rho_i}(x_i)})=|x-x_i|-\rho_i\ge\frac12|x-x_i|\ge\frac12 d(x)$.
	For $x\in F$ and $i\notin\mathcal I'$, \eqref{e:boundary balls} gives $\dist(x,\overline{B_{\rho_i}(x_i)})\ge\frac18-2\delta_0>\frac1{10}$.
	Therefore
	\[
	\dist(x,\partial U)\ge\min\{\tfrac1{10},\tfrac12d(x)\},\qquad d(x)\ge r_0\qquad\text{for all }x\in F .
	\]
	Using $|x|\le2$ and $|\nabla\varphi|\le36$, we get
	\[
	\int_U\frac{|x|^2|\nabla\varphi|^2}{\dist^2(x,\partial U)}\,dx\le C\int_F\Big(1+\frac1{d^2(x)}\Big)dx
	\le C+C\int_{F\cap\{d\ge r_{K-1}\}}\frac{dx}{d^2}+C\sum_{k=1}^{K-1}\int_{F\cap U_k}\frac{dx}{d^2},
	\]
	where $U_k:=\{x:\ r_{k-1}\le d(x)<r_k\}$.
	The second term is bounded by $C|B_2|\,r_{K-1}^{-2}\le C M^4c_*^{-2}$.
	
	For the third term we claim that, for $1\le k\le K-1$,
	\begin{equation}\label{e:covering}
		F\cap U_k\subseteq\bigcup_{z\in D_{k+1}}B_{2r_{k+1}}(z).
	\end{equation}
	Let $x\in F\cap U_k$, let $x_i\in I$ be a closest center to $x$, so that $|x-x_i|=d(x)<r_k$, and let $z:=y_{k+1}(x_i)\in I_{k+1}$.
	
	If $z\notin D_{k+1}$, then $\ch(z)=\{z\}$ and hence $y_k(x_i)=z$.
	By \eqref{e:ancestors}, $|x-z|<r_k+2r_k<4r_k$, i.e.\ $x\in B_{2s}(z)$ with $(z,s)=(z,2r_k)\in\mathcal A$.
	By (ii) above, $x\in G$, which contradicts $x\in F$.
	Hence $z\in D_{k+1}$, and by \eqref{e:ancestors}
	\[
	|x-z|<r_k+\tfrac{M}{M-1}r_{k+1}<2r_{k+1}.
	\]
	This proves \eqref{e:covering}.
	
	Since $d\ge r_{k-1}$ on $U_k$, \eqref{e:covering} and \eqref{e:telescopic} give
	\[
	\sum_{k=1}^{K-1}\int_{F\cap U_k}\frac{dx}{d^2}
	\le\sum_{k=1}^{K-1}\frac{4\pi r_{k+1}^2}{r_{k-1}^2}\,\#D_{k+1}
	=4\pi M^4\sum_{k=1}^{K-1}\#D_{k+1}\le4\pi M^4N' .
	\]
	Collecting the estimates proves \eqref{e:foliation estimate}.
	
	If $N'=0$, we take $\varphi:=\varphi_0$ as above. Then (1) and (2) are clear. By \eqref{e:boundary balls}, every ball lies where $\varphi_0\in\{0,1\}$ is locally constant, which gives (3). The gradient $\nabla\varphi_0$ is supported in $\{\frac54\le|x|\le\frac74\}$, where $\dist(x,\partial U)\ge\frac14-2\delta_0$, which gives (4).	
\end{proof}

\subsection{Energy transfer} In the following we use the notation $\A_r(x_0) := B_{2r}\setminus \overline{B_r}(x_0)$ for the annuli. 


\begin{lemma}[Energy transfer]\label{l:energy transfer}
	There exists a universal constant $C>0$ with the following property.
	Let $\mu\in[\frac12,1)$.
	Let $R>0$, $x_0\in\mathbb R^2$ and $A\in L^2(B_{2R}(x_0);\mathbb R^{2\times2})$ with $\curl A\in L^\infty(B_{2R}(x_0))$.
	Let $B_i=B_{\rho_i}(x_i)$, $i=1,\dots,N$, be open balls with pairwise disjoint closures such that
	\begin{itemize}
		\item[(a)] $\operatorname{spt}(\curl A)\cap B_{2R}(x_0)\subset\bigcup_{i}B_i$;
		\item[(b)] $\partial B_{2R}(x_0)\cap\overline{B_i}=\emptyset$ for every $i$;
		\item[(c)] $\sum\big\{\rho_i : B_i\cap \A_{R}
		(x_0)\neq\emptyset\big\}<\delta_0R$, with $\delta_0$ the constant of Lemma~\ref{l:foliation}.
	\end{itemize}
	Assume moreover that $A$ is harmonic in $B_{2R}(x_0)\setminus\bigcup_i\overline{B_i}$.
	Set
	\[
	I':=\{i : B_i\subset \A_R(x_0)\},\qquad N':=\# I',\qquad U:=\A_R(x_0)\setminus\bigcup_{i}\overline{B_i}.
	\]
	Then, for every $\eta>0$,
	\begin{equation}\label{e:energy transfer}
		\big|\curl A\big(B_{2R}(x_0)\big)\big|
		\le\sum_{i\in I'}\big|\curl A(B_i)\big|
		+\frac{C}{\eta}\int_{U}\operatorname{dist}^2_{K_\mu}(A)\,\d x+\eta\,(1+N').
	\end{equation}
\end{lemma}

\begin{proof}
	Without loss of generality we assume $x_0=0$ and let $\varphi$ be the function given by Lemma~\ref{l:foliation} for the balls $\{B_i\}$. For $i\in I'$ we denote by $c_i\in[0,1]$ the constant value of $\varphi$ on $B_i$.
	Since $\varphi=1$ on $B_{9R/8}$, $\varphi=0$ on $B_{2R}\setminus\overline{B_{15R/8}}$ and $\varphi$ is constant on a neighborhood of each $\overline{B_i}$, the support of $\nabla\varphi$ is a compact subset of $U$, where $A$ is harmonic and hence smooth.
	
	We first show that
	\begin{equation}\label{e:ibp}
		\int_{B_{2R}}\varphi\,\curl A\, \d x=\int_U\big(\nabla A\cdot x\big)\,J\nabla\varphi\,\d x ,
	\end{equation}
	where $(\nabla A\cdot v)_{ik}:=\sum_j\partial_ka_{ij}v_j$.
	By definition of the distributional curl, $$\int\psi\,\curl A\,dx=-\int A\,J\nabla\psi\,dx \qquad \forall\,\psi\in C^\infty_c(B_{2R}),$$ and since $A\in L^2$ and $\curl A\in L^\infty$ the identity extends by mollification to every Lipschitz $\psi$ with compact support in $B_{2R}$, in particular to $\psi=\varphi$.
	On the other hand, $J\nabla\varphi$ is divergence-free and has compact support in $U$, so that, integrating by parts the smooth field $A(x) x$, whose gradient is $A+\nabla A\cdot x$,
	\[
	0=\int_U\nabla\big[A(x) x\big]\,J\nabla\varphi\,dx=\int A\,J\nabla\varphi\,\d x+\int_U\big(\nabla A\cdot x\big)\,J\nabla\varphi\,\d x ,
	\]
	and \eqref{e:ibp} follows.
	
	Next we compute the left-hand side of \eqref{e:ibp}. By (b), each $B_i$ is either contained in $B_{2R}$ or has closure disjoint from $\overline{B_{2R}}$. If $B_i\subset B_{2R}$ and $i\notin I'$, then 
	by $\rho_i<\delta_0R$ we get $B_i\subset B_{(1+2\delta_0)R}\subset B_{9R/8}$; hence, $\varphi=1$ on $B_i$. By (a) and the disjointness of the balls,
	\[
	\int_{B_{2R}}\varphi\,\curl A\,dx
	=\curl A(B_{2R})-\sum_{i\in I'}(1-c_i)\,\curl A(B_i),
	\]
	which, together with \eqref{e:ibp} and $c_i\in[0,1]$, gives
	\begin{equation}\label{e:step3}
		\big|\curl A(B_{2R})\big|\le\sum_{i\in I'}\big|\curl A(B_i)\big|
		+\int_U|\nabla A|\,|x|\,|\nabla\varphi|\,\d x .
	\end{equation}
	Let $d(x):=\operatorname{dist}(x,\partial U)$. By Young's inequality and \eqref{e:foliation estimate}, for every $\eta>0$
	\begin{align}\label{e:young}
		\int_U|\nabla A|\,|x|\,|\nabla\varphi|\,\d x&
		\le\frac1{4\eta}\int_U|\nabla A|^2d^2\,dx+\eta\int_U\frac{|x|^2|\nabla\varphi|^2}{d^2}\,\d x \notag\\
		& \le\frac1{4\eta}\int_U|\nabla A|^2d^2\,\d x+C\eta\,(1+N').
	\end{align}
	
	It remains to estimate the first term on the right-hand side by the elastic energy. Consider a Whitney-type covering of $U$, i.e.\ a countable family of balls $B_{r_k}(y_k)$ with $4r_k=d(y_k)$, such that the balls $B_{r_k/2}(y_k)$ cover $U$ and every point of $U$ belongs to at most $M$ of the balls $B_{2r_k}(y_k)$, with $M$ universal. Then $B_{2r_k}(y_k)\subset U$ and $d\le\frac92r_k$ on $B_{r_k/2}(y_k)$.
	On $B_{2r_k}(y_k)$ the field $A$ is harmonic and curl-free, hence the gradient of a $W^{1,2}$ function, and Theorem~\ref{t:uniform rigidity} gives $S_k\in K_\mu$ with
	\[
	\int_{B_{3r_k/4}(y_k)}|A-S_k|^2\,\d x\le C\int_{B_{2r_k}(y_k)}\dist_{K_\mu}^2(A)\,\d x ,
	\]
	with $C$ universal, in particular independent of $k$ and of $\mu\in[\tfrac12,1)$. Since $A-S_k$ is harmonic, the Caccioppoli inequality gives
	\[
	\int_{B_{r_k/2}(y_k)}|\nabla A|^2d^2\,\d x
	\le\frac{81}4\,r_k^2\cdot\frac{64}{r_k^2}\int_{B_{3r_k/4}(y_k)}|A-S_k|^2\,\d x
	\le C\int_{B_{2r_k}(y_k)}\dist^2_{K_\mu}(A)\,\d x ,
	\]
	and summing over $k$ we obtain
	\begin{equation}\label{e:whitney}
		\int_U|\nabla A|^2d^2\,\d x\le C\,M\int_U\dist_{K_\mu}^2(A)\,\d x .
	\end{equation}
	The conclusion \eqref{e:energy transfer} follows from \eqref{e:step3}, \eqref{e:young} and \eqref{e:whitney}, after replacing $\eta$ with $\eta/C$.
\end{proof}

\begin{remark}\label{r:energy transfer mu}
	The restriction $\mu\in[\tfrac12,1)$ is used only to have a rigidity constant independent of $\mu$. For a fixed $\mu\in(0,1)$, the same proof with Corollary~\ref{c:rigidity} in place of Theorem~\ref{t:uniform rigidity} gives \eqref{e:energy transfer} with a constant depending on $\mu$: for curl-free fields both sides of \eqref{e:rigidity0} scale in the same way, so the rigidity constant on every ball equals the one on $B_1$.
\end{remark}

\section{Logarithmic estimate}\label{s:log}

In this section we give an abstract formulation of the logarithmic estimate of
Lauteri--Luckhaus~\cite{LauteriLuckhaus2016}. The only ingredient coming from the elastic model is
hypothesis (H2) below, which is verified in Section~6 by means of Lemma~\ref{l:energy transfer}.

\subsection{Notation and two covering lemmas}
For an open ball $B$ we write $x_B$ and $\rho_B$ for its center and radius. 
We call a finite family $\F$ of open balls with pairwise disjoint closures a \emph{disjoint balls family} and we write
\[
\cup\F:=\bigcup_{B\in\F}B,\qquad
\partial\F:=\bigcup_{B\in\F}\partial B,\qquad
\Sigma(\F):=\sum_{B\in\F}\rho_B,\qquad
\#\F := \text{cardinality of $\F$}.
\]
A function $\nu\in L^1_{\rm loc}$ is identified
with the measure $\nu\,\Le$; in particular $\nu(\partial E)=0$ for every Lipschitz
set $E$. Finally, $\delta_0=1/4000$ is the absolute constant of
Lemma~\ref{l:foliation}.

\begin{lemma}[Merging]\label{l:merge}
	Let $\mathcal B$ be a finite family of open balls in $\R^2$. Then there exists a disjoint balls family $\mathcal M$ (i.e. with disjoint closures), called \emph{the merging} of $\mathcal B$, such that
	\begin{enumerate}[label=(\roman*)]
		\item every ball of $\mathcal B$ is contained in some ball of $\mathcal M$;
		\item $\Sigma(\mathcal M)\le \Sigma(\mathcal B)$ and $\#\mathcal M\le \#\mathcal B$.
	\end{enumerate}
\end{lemma}

\begin{proof}
	If two balls $B_r(x),B_s(y)$ of the current family have intersecting closures,
	i.e.\ $|x-y|\le r+s$, we replace them with the smallest open ball $B_t(z)$
	containing both: if one of them contains the other, $B_t(z)$ is the larger one;
	otherwise
	\[
	z=\frac{x+y}{2}+\frac{(s-r)(y-x)}{2|y-x|},\qquad t=\frac{|x-y|+r+s}{2}\le r+s .
	\]
	Each step decreases the cardinality by one, does not increase the sum of the
	radii and preserves (i). The procedure stops after at most $\#\mathcal B-1$
	steps, when the closures are pairwise disjoint.
\end{proof}

The following is a variant of the Besicovitch-type lemma of \cite{LauteriLuckhaus2016}, in which the selected balls are
disjoint after a dilation by $3$ (instead of $2$). This extra room is needed to
choose \emph{good radii} in the proof of Theorem~\ref{t:log}. The proof is given in
Appendix~B.

\begin{lemma}[Besicovitch-type selection]\label{l:besi}
	Let $B_R(y)\subset\R^n$, let $\{B_{\rho_i}(x_i)\}_{i\in I}$ be a finite family of
	balls with $B_{16\rho_i}(x_i)\subset B_R(y)$ for all $i\in I$, and let $\zeta$ be a
	finite positive Borel measure with
	$\zeta\big(B_R(y)\setminus\bigcup_{i}B_{\rho_i}(x_i)\big)=0$. Then there exist
	$\tilde I\subset I$ and radii $R_i\ge 2\rho_i$, $i\in\tilde I$, such that
	\begin{enumerate}[label=(\arabic*)]
		\item the balls $B_{3R_i}(x_i)$, $i\in\tilde I$, are pairwise disjoint and contained in $B_R(y)$;
		\item $\displaystyle\sum_{i\in\tilde I}\zeta\big(B_{R_i}(x_i)\big)\ge c_n\,\zeta\big(B_R(y)\big)$, with $c_n:=107^{-n}$.
	\end{enumerate}
\end{lemma}

\subsection{The abstract estimate}

\begin{theorem}\label{t:log}
	Let $\lambda,\tau_0,C_0>0$. There exist $\delta_1\in(0,1)$ depending only on
	$\lambda$, and $C>0$ depending on $\lambda,\tau_0,C_0$, with the following property.
	
	Let $\eps>0$, $\nu\in L^\infty(\R^2;\R^2)$, let $S\subset\R^2$ be a closed set with $\supp\nu\subset S$ and let $\sigma$ be a positive Radon measure on $\R^2$. Set
	\[
	\mu_0^\eps:=\eps^{-1}\,\Le\res B_{\lambda\eps}(S),\qquad
	\mu_1^\eps:=\eps^{-1}\sigma .
	\]
	Assume:
	\begin{enumerate}[label=\textup{(H\arabic*)},start=0,leftmargin=*]
		\item \textbf{$L^\infty$ bound:} $\|\nu\|_\infty\le\eps^{-1}$;
		\item \textbf{Quantization:} for every bounded open set $E\subset\R^2$ whose boundary is a closed Lipschitz curve with $\partial E\cap S=\emptyset$,
		\[
		\text{either}\quad \nu(E)=0 \quad \text{or}\quad  |\nu(E)|\ge\tau_0\eps;
		\]
		\item \textbf{Energy transfer:} for every $x_0\in\R^2$, $\rho>0$ and every disjoint balls family $\F$ such that
		\begin{enumerate}[label=(\alph*)]
			\item $S\cap\ov B_{2\rho}(x_0)\subset\cup\F$,
			\item $\partial B_{2\rho}(x_0)\cap\ov B=\emptyset$ for every $B\in\F$,
			\item $\sum\big\{\rho_B:\ B\in\F,\ \ov B\cap\ov \A_\rho(x_0)\neq\emptyset\big\}<\delta_0\rho$,
		\end{enumerate}
		we have, for every $\eta>0$,
		\begin{equation}\label{e:H2}
			|\nu(B_{2\rho}(x_0))|\le\sum_{B\in\F_A}|\nu(B)|
			+\frac{C_0}{\eta}\,\sigma\big(\A_\rho(x_0)\setminus\cup\F\big)+\eta\,(1+\#\F_A),
		\end{equation}
		where $\F_A:=\{B\in\F:\ B\subset \A_\rho(x_0)\}$.
	\end{enumerate}
	Then, for every $x_0\in\R^2$ and every $R\ge 8\eps$, there exists a Borel set
	$E\subset(R/2,5R/8)$ with $|E|\ge R/16$ such that
	\begin{equation}\label{e:main}
		|\nu(B_r(x_0))|\le C\,\omega\!\left(\frac{\mu_0^\eps(B_R(x_0))}{R}\right)
		\big(\mu^\eps_0+\mu^\eps_1\big)(B_R(x_0))\qquad\forall\,r\in E,
	\end{equation}
	where $\omega(t):=|\log\min\{t,\delta_1\}|^{-1}$ for $t>0$ and $\omega(0):=0$.
	If moreover $\mu_0^\eps(B_R(x_0))\le\delta_1 R$, then $E$ can be chosen so that,
	in addition, $\partial B_r(x_0)\cap S=\emptyset$ for every $r\in E$.
\end{theorem}

\begin{remark}[Scaling]\label{r:scaling}
	Given $x_0\in\R^2$ and $R>0$, set
	\[
	\nu'(x):=R\,\nu(x_0+Rx),\qquad S':=R^{-1}(S-x_0),\qquad \sigma'(E):=R^{-2}\sigma(x_0+RE),\qquad \eps':=\eps/R .
	\]
	Then $\|\nu'\|_\infty\le (\eps')^{-1}$, and 
	\[
	\nu'(E) = R^{-1}\nu (x_0+RE), \quad \mu_0^{\eps'}{}'(E)=R^{-1}\mu_0^\eps(x_0+RE),\quad  \mu_1^{\eps'}{}'(E)=R^{-1}\mu_1^\eps(x_0+RE).
	\]
	Therefore, (H0)--(H2) hold for $(\nu',S',\sigma',\eps')$ with the same constants
	(in (H2) one replaces $\eta$ with $\eta/R$). Both sides of \eqref{e:main} scale by
	the factor $R$. Hence it suffices to prove Theorem~\ref{t:log} for $x_0=0$ and $R=1$,
	i.e.\ for $\eps\le 1/8$.
\end{remark}

\subsection{Sparse scales and multiscale coverings}

\begin{definition}\label{d:sparse}
	Let $\F$ be a family of balls and $B\in\F$. For $t>0$ we set
	\[
	\ell_\F(B,t):=\Hu\big(\partial\F\cap\ov \A_t(x_B)\big),
	\qquad
	\rb_\F(B):=\sup\Big(\{\rho_B\}\cup\{t>0:\ \ell_\F(B,t)>\delta_0 t\}\Big).
	\]
\end{definition}

Since $\ell_\F(B,t)\le\Hu(\partial\F)=2\pi\Sigma(\F)$, we have
$\rb_\F(B)\le\max\{\rho_B,2\pi\Sigma(\F)/\delta_0\}<\infty$.

\begin{lemma}[Sparse scales]\label{l:sparse}
	Let $\F$ be a disjoint balls family, $B\in\F$ and $\rb:=\rb_\F(B)$.
	\begin{enumerate}[label=(\alph*)]
		\item $\rho_B\le\rb\le\delta_0^{-1}\,\Hu\big(\partial\F\cap\Sigma_\F(B)\big)$, where
		$\Sigma_\F(B):=\ov B_{2\rb}(x_B)\setminus B_{\rb/2}(x_B)$.
		\item For every $\rho>2\rb$, every $B'\in\F$ with $\ov{B'}\cap\ov \A_\rho(x_B)\neq\emptyset$
		is different from $B$, and
		\[
		\sum\big\{\rho_{B'}:\ B'\in\F,\ \ov{B'}\cap\ov \A_\rho(x_B)\neq\emptyset\big\}
		\le\frac{7}{4\pi}\,\delta_0\rho<\delta_0\rho .
		\]
	\end{enumerate}
\end{lemma}

\begin{proof}
	(a) The first inequality holds by definition. If $\rb=\rho_B$, then
	$\partial B\subset\Sigma_\F(B)$, so that
	$\Hu(\partial\F\cap\Sigma_\F(B))\ge2\pi\rho_B\ge\delta_0\rb$. If $\rb>\rho_B$,
	there are $t_n\le\rb$ with $t_n\to\rb$ and $\ell_\F(B,t_n)>\delta_0t_n$; for $n$
	large $t_n\ge\rb/2$, hence $\ov \A_{t_n}(x_B)\subset\Sigma_\F(B)$ and
	$\delta_0t_n<\Hu(\partial\F\cap\Sigma_\F(B))$. Letting $n\to\infty$ gives the claim.
	
	(b) Let $\rho>2\rb$ and $W:=\ov B_{4\rho}(x_B)\setminus B_{\rho/2}(x_B)
	=\ov \A_{\rho/2}(x_B)\cup\ov \A_{\rho}(x_B)\cup\ov \A_{2\rho}(x_B)$. Since $\rho/2>\rb$,
	by definition of $\rb$
	\begin{equation}\label{e:W}
		\Hu(\partial\F\cap W)\le\delta_0\Big(\frac\rho2+\rho+2\rho\Big)=\frac72\,\delta_0\rho .
	\end{equation}
	Let $B'\in\F$ with $\ov{B'}\cap\ov \A_\rho(x_B)\neq\emptyset$ and write $\rho':=\rho_{B'}$.
	Since $\ov B=\ov B_{\rho_B}(x_B)$ and $\rho_B\le\rb<\rho$, $B'\neq B$; as the
	closures are disjoint and $x_B\in B$, we get $x_B\notin\ov{B'}$. We claim that
	\begin{equation}\label{e:claim}
		\Hu(\partial B'\cap W)\ge\min\{2\pi\rho',\rho/4\}.
	\end{equation}
	Fix $p\in\ov{B'}\cap\ov \A_\rho(x_B)$.
	
	\emph{Case $\rho'\le\rho/8$.} Then $\ov{B'}\subset\ov B_{\rho/4}(p)\subset
	\{3\rho/4\le|y-x_B|\le 9\rho/4\}\subset W$, and $\Hu(\partial B'\cap W)=2\pi\rho'$.
	
	\emph{Case $\rho'>\rho/8$ and $\partial B'\cap B_{\rho/2}(x_B)\neq\emptyset$.}
	The point of $\ov{B'}$ farthest from $x_B$ lies on $\partial B'$, at distance
	$\ge|p-x_B|\ge\rho$. Since $\partial B'$ is connected, it contains an arc joining
	$\partial B_{\rho/2}(x_B)$ to $\partial B_\rho(x_B)$ inside
	$\{\rho/2\le|y-x_B|\le\rho\}\subset W$, whose length is at least $\rho/2$.
	
	\emph{Case $\rho'>\rho/8$ and $\partial B'\cap B_{\rho/2}(x_B)=\emptyset$.}
	The connected set $B_{\rho/2}(x_B)$ does not meet $\partial B'$ and contains
	$x_B\notin\ov{B'}$, hence it is disjoint from $\ov{B'}$. Let $q\in\ov{B'}$ be a
	point closest to $x_B$: then $q\in\partial B'$, $\rho/2\le|q-x_B|\le|p-x_B|\le2\rho$,
	and $|y-x_B|\ge\rho/2$ for all $y\in\ov{B'}$. A circle of radius $r\ge s$ through
	$q$ has length $4r\arcsin\frac{s}{2r}\ge 2s$ inside $B_s(q)$; with $s=\rho/8$ we get
	that $\partial B'\cap B_{\rho/8}(q)$ has length $\ge\rho/4$ and is contained in
	$\{\rho/2\le|y-x_B|\le 17\rho/8\}\subset W$. This proves \eqref{e:claim}.
	
	By \eqref{e:W} and \eqref{e:claim}, $\min\{2\pi\rho',\rho/4\}\le\frac72\delta_0\rho<\rho/4$
	because $\delta_0<1/14$; hence $\Hu(\partial B'\cap W)\ge2\pi\rho'$ for every
	such $B'$. Since the circles $\partial B'$ are pairwise disjoint, summing and
	using \eqref{e:W} again gives $2\pi\sum\rho'\le\frac72\delta_0\rho$.
\end{proof}

\begin{proposition}[One step of the multiscale covering]\label{p:step}
	Let $\F$ be a family of balls and let $\rb(B):=\rb_\F(B)$ for $B\in\F$. Let
	$\hat\F\subset\F$ be obtained greedily: order $\F$ by nonincreasing $\rb$, and
	select $B$ if $\ov B_{2\rb(B)}(x_B)$ is disjoint from $\ov B_{2\rb(\hat B)}(x_{\hat B})$
	for all previously selected $\hat B$. Let $\F'$ be the merging of
	$\{B_{36\rb(\hat B)}(x_{\hat B}):\ \hat B\in\hat\F\}$. Then:
	\begin{enumerate}[label=(\roman*)]
		\item for every $B\in\F$ there is exactly one $P\in\F'$ with $B\subset P$ (called the
		\emph{parent} of $B$; $B$ is called a \emph{child} of $P$); moreover
		$B_{32\rb(B)}(x_B)\subset P$. Every $P\in\F'$ has at least one child and $\cup\F\subset \cup\F'$;
		\item $\Sigma(\F')\le\Lambda\,\Sigma(\F)$, where $\Lambda:=72\pi/\delta_0$.
	\end{enumerate}
\end{proposition}

\begin{proof}
	(i) If $B\in\F$ is not selected, $\ov B_{2\rb(B)}(x_B)$ meets
	$\ov B_{2\rb(\hat B)}(x_{\hat B})$ for some $\hat B\in\hat\F$ with
	$\rb(\hat B)\ge\rb(B)$; thus $|x_B-x_{\hat B}|\le4\rb(\hat B)$ and
	$B_{32\rb(B)}(x_B)\subset B_{36\rb(\hat B)}(x_{\hat B})$. The same inclusion is
	trivial if $B=\hat B$ is selected. By Lemma~\ref{l:merge}(i),
	$B_{36\rb(\hat B)}(x_{\hat B})$ is contained in some $P\in\F'$, and
	$B\subset B_{32\rb(B)}(x_B)\subset P$ because $\rho_B\le\rb(B)$. Uniqueness of $P$
	follows from the disjointness of the closures of the balls of $\F'$. Each
	$P\in\F'$ contains some $B_{36\rb(\hat B)}(x_{\hat B})\supset\hat B$.
	
	(ii) The sets $\Sigma_\F(\hat B)\subset\ov B_{2\rb(\hat B)}(x_{\hat B})$,
	$\hat B\in\hat\F$, are pairwise disjoint. By Lemma~\ref{l:merge}(ii) and
	Lemma~\ref{l:sparse}(a),
	\[
	\Sigma(\F')\le36\sum_{\hat B\in\hat\F}\rb(\hat B)
	\le\frac{36}{\delta_0}\sum_{\hat B\in\hat\F}\Hu\big(\partial\F\cap\Sigma_\F(\hat B)\big)
	\le\frac{36}{\delta_0}\Hu(\partial\F)=\frac{72\pi}{\delta_0}\,\Sigma(\F).\qedhere
	\]
%
\end{proof}

\subsection{Proof of Theorem~\ref{t:log}}

By Remark~\ref{r:scaling} we may assume $x_0=0$, $R=1$, $\eps\le1/8$. We write
$\mu_0:=\mu_0^\eps$, $\mu_1:=\mu_1^\eps$, and fix
\[
\theta:=\frac1{32},\qquad \ell:=\min\{\lambda,1\},\qquad
\Lambda:=\frac{72\pi}{\delta_0},\qquad \delta_1:=\Big(\frac{\ell\theta}{\Lambda}\Big)^4 .
\]

We first discuss two trivial cases. If $\mu_0(B_1)=0$, then $S\cap B_1=\emptyset$, so $\nu=0$ a.e.\ on $B_1$ and $E:=(1/2,5/8)$ works. 
If $\mu_0(B_1)>\delta_1$, then by (H0) $|\nu(B_r)|\le|\nu|(B_1)\le \mu_0(B_1)=|\log\delta_1|\,\omega(\mu_0(B_1))\,\mu_0(B_1)$ for every $r$, and
again $E:=(1/2,5/8)$ works provided $C\ge|\log\delta_1|$. From now on
$0<\mu_0(B_1)\le\delta_1$.

Let $Z$ be a maximal subset of the compact set $S\cap\ov B_{7/8}$ with
$|z-z'|\ge2\ell\eps$ for distinct $z,z'\in Z$. The balls $B_{\ell\eps}(z)$, $z\in Z$,
are pairwise disjoint and contained in $B_{\lambda\eps}(S)\cap B_1$ (since
$\ell\le\lambda$ and $\ell\eps\le1/8$), so
\[
\#Z\,\pi\ell^2\eps^2\le|B_{\lambda\eps}(S)\cap B_1|=\eps \mu_0(B_1),
\qquad\text{i.e.}\qquad \#Z\le\frac{\mu_0(B_1)}{\pi\ell^2\eps}.
\]
By maximality, $S\cap\ov B_{7/8}\subset\bigcup_{z\in Z}B_{2\ell\eps}(z)$. Let $\F_0$
be the merging of $\{B_{2\ell\eps}(z)\}_{z\in Z}$. By Lemma~\ref{l:merge},
\begin{equation}\label{e:F0}
	S\cap\ov B_{7/8}\subset\cup \F_0,\qquad
	\#\F_0\le\frac{\mu_0(B_1)}{\pi\ell^2\eps},\qquad
	\Sigma(\F_0)\le2\ell\eps\,\#Z\le\frac{\mu_0(B_1)}{\ell}.
\end{equation}

Let $\F_{k+1}:=\F_k'$ be obtained from $\F_k$ as in Proposition~\ref{p:step},
and let
\[
K:=\max\big\{k\ge0:\ \Lambda^k \mu_0(B_1)/\ell\le\theta\big\},
\]
which is well defined since $\mu_0(B_1)/\ell\le\delta_1/\ell\le\theta$. By
Proposition~\ref{p:step}(ii) and \eqref{e:F0}, $\Sigma(\F_k)\le\theta$ for
$k\le K$; in particular every ball of $\F_k$, $k\le K$, has radius at most
$\theta$. Moreover, since $\mu_0(B_1)\le\delta_1$ gives $|\log(\ell\theta)|\le\frac14|\log \mu_0(B_1)|$
and $\log\Lambda\le\frac14|\log \mu_0(B_1)|$,
\begin{equation}\label{e:K}
	K\ge\frac{\log(\ell\theta/\mu_0(B_1))}{\log\Lambda}-1
	\ge\frac{\frac34|\log \mu_0(B_1)|}{\log\Lambda}-\frac{\frac14|\log \mu_0(B_1)|}{\log\Lambda}
	=\frac{|\log \mu_0(B_1)|}{2\log\Lambda} \ge2.
\end{equation}
By \eqref{e:F0} and Proposition~\ref{p:step}(i),
$S\cap\ov B_{7/8}\subset\cup \F_k$ for every $k$. We set
\[
\G_k:=\{B\in\F_k:\ B\subset B_{3/4}\},\qquad
\tau_k:=\sum_{B\in\G_k}\frac{|\nu(B)|}{|B|}\,\Le\res B,\qquad
\|\tau_k\|:=\tau_k(\R^2)=\sum_{B\in\G_k}|\nu(B)|,
\]
and denote by $\operatorname{ch}(B)$ the set of children of $B\in\F_{k+1}$.
Let $B\in\G_{k+1}$. 
By construction all children of $B$ are contained in $B\subset B_{3/4}$ and belong to $\G_k$, and, since $S\cap B\subset S\cap\ov B_{7/8}\subset\cup\F_k$, we have $\nu(B)=\sum_{B'\in\operatorname{ch}(B)}\nu(B')$. Consequently

\begin{enumerate}[label=(O\arabic*),leftmargin=*] 
	\item\label{o1} $\|\tau_{k+1}\|\le\|\tau_k\|$,
	\item\label{o2} $\sum_{B\in\G_{k+1}}\tau_k(B)\le\|\tau_k\|$,
	\item\label{o3} $\|\tau_0\|\le|\nu|(B_{3/4})\le \mu_0(B_1)$, by (H0).
\end{enumerate}

The key estimate is the following one-step decay: there exist constants $q\in (0,1)$ and $C_1>0$ with the following property: for every $0\le k<K$ and every $B\in\G_{k+1}$, there exist a Borel set
$D_B\subset B\setminus\cup \F_k$ and an integer $N_B$ with
$N_B\le2(\#\operatorname{ch}(B)-1)$ such that
\begin{equation}\label{e:onestep}
	|\nu(B)|\le q\,\tau_k(B)+C_1\,\mu_1(D_B)+C_1\,\eps N_B .
\end{equation}

To prove it, let $q:= 1- c_2/2$ with $c_2:=107^{-2}$ the constant of Lemma~\ref{l:besi} for $n=2$, and $C_1:=\frac{4C_0}{c_2\tau_0}+\tau_0$.
Let $B_i$, $i\in I_B$, be the children of $B$, with centers $x_i$, radii $\rho_i$,
and $\rb_i:=\rb_{\F_k}(B_i)$. By Proposition~\ref{p:step}(i),
$B_i\subset B_{2\rb_i}(x_i)\subset B_{32\rb_i}(x_i)\subset B$.
By~\ref{o1}, $\tau_k\res B$ is concentrated on $\bigcup_iB_i$. Lemma~\ref{l:besi}, applied in the
ball $B$ to the balls $B_{2\rb_i}(x_i)$ and to $\zeta:=\tau_k\res B$, gives
$I'_B\subset I_B$ and radii $R_i\ge4\rb_i$ such that
\begin{gather}
	B_{3R_i}(x_i)\;\text{ are pairwise disjoint for $i\in I'_B$ and are contained in }B,\label{e:besi-1}\\
	\sum_{i\in I'_B}\tau_k\big(B_{R_i}(x_i)\big)\ge c_2\,\tau_k(B).\label{e:besi-2}
\end{gather}

We now choose good radii. Fix $i\in I'_B$. If $r\in[R_i,3R_i/2]$, every
$B'\in\F_k$ whose closure meets $\partial B_{2r}(x_i)$ meets
$\ov \A_{2R_i}(x_i)$; since $2R_i>2\rb_i$, Lemma~\ref{l:sparse}(b) gives that the
sum of the radii of these balls is less than $2\delta_0R_i$. For each such $B'$,
the set of $r$ with $\partial B_{2r}(x_i)\cap\ov{B'}\neq\emptyset$ is an interval
of length at most $2\rho_{B'}$. Hence the set of bad $r\in[R_i,3R_i/2]$ has measure
less than $4\delta_0R_i<R_i/2$, and we can fix $r_i\in[R_i,3R_i/2]$ such that:
\begin{enumerate}[label=($\alpha$\arabic*),leftmargin=*]
	\item $\ov B_{2r_i}(x_i)\subset B_{3R_i}(x_i)\subset B$, and these balls are
	pairwise disjoint for $i\in I_B'$;
	\item $\partial B_{2r_i}(x_i)\cap\ov{B'}=\emptyset$ for all $B'\in\F_k$; since
	$S\cap B\subset\cup\F_k$, also $\partial B_{2r_i}(x_i)\cap S=\emptyset$;
	\item $r_i\ge R_i>2\rb_i$.
\end{enumerate}

We can then apply (H2). Let $\A_i:=\A_{r_i}(x_i)$ and
$\F^{(i)}:=\{B'\in\F_k:\ \ov{B'}\cap\ov B_{2r_i}(x_i)\neq\emptyset\}$, which by
($\alpha$1) and~\ref{o1} consists of children of $B$. Condition (a) of (H2)
holds because $S\cap\ov B_{2r_i}(x_i)\subset S\cap B\subset\cup \F_k$, (b) is
($\alpha$2), and (c) follows from ($\alpha$3) and Lemma~\ref{l:sparse}(b). Let
$N_i:=\#\{B'\in\F_k:\ B'\subset \A_i\}$. Since these balls are children of $B$,
they belong to $\G_k$ and $\sum_{B'\subset \A_i}|\nu(B')|\le\tau_k(\A_i)$. Moreover
$\A_i\cap\cup\F_k=\A_i\cap\cup \F^{(i)}$. Hence (H2) with $\eta:=\frac{c_2}{2}\tau_0\eps$
gives
\[
|\nu(B_{2r_i}(x_i))|\le\tau_k(\A_i)+\frac{C_0}{\frac{c_2}{2}\tau_0}\,
\mu_1\big(\A_i\setminus\cup \F_k\big)+\frac{c_2}{2}\tau_0\eps\,(1+N_i).
\]
By ($\alpha$2) and (H1), either $\nu(B_{2r_i}(x_i))=0$ or
$\frac{c_2}{2}\tau_0\eps\le\frac{c_2}{2}|\nu(B_{2r_i}(x_i))|$. In both cases
\begin{equation}\label{e:absorbed}
	q\,|\nu(B_{2r_i}(x_i))|\le\tau_k(\A_i)+\frac{C_0}{\frac{c_2}{2}\tau_0}\,
	\mu_1\big(\A_i\setminus\cup \F_k\big)+\frac{c_2}{2}\tau_0\eps N_i .
\end{equation}

By ($\alpha$2), every child of $B$ is either contained in some
$B_{2r_i}(x_i)$ or disjoint from all of them. By~\ref{o1},
\[
|\nu(B)|\le\sum_{i\in I'_B}|\nu(B_{2r_i}(x_i))|+\tau_k\Big(B\setminus\bigcup_{i\in I_B'}B_{2r_i}(x_i)\Big).
\]
Since $\tau_k\ll\Le$ and $r_i\ge R_i$,
$\tau_k(\A_i)=\tau_k(B_{2r_i}(x_i))-\tau_k(B_{r_i}(x_i))\le\tau_k(B_{2r_i}(x_i))-\tau_k(B_{R_i}(x_i))$.
Plugging \eqref{e:absorbed} in, and using ($\alpha$1) and \eqref{e:besi-2},
\begin{align*}
q|\nu(B)|&\le\tau_k(B)-\sum_{i\in I'_B}\tau_k(B_{R_i}(x_i))
+\frac{C_0}{\frac{c_2}{2}\tau_0}\mu_1(D_B)+\frac{c_2}{2}\tau_0\eps N_B\\
&\le(1-c_2)\tau_k(B)+\frac{C_0}{\frac{c_2}{2}\tau_0}\mu_1(D_B)+\frac{c_2}{2}\tau_0\eps N_B,
\end{align*}
where $D_B:=\bigcup_{i\in I'_B}(\A_i\setminus\cup\F_k)\subset B\setminus\cup \F_k$ and $N_B:=\sum_{i\in I'_B}N_i$. Since $\frac{1-c_2}{1-c_2/2}\le1-\frac{c_2}{2}=q$
and $\frac1{1-\frac{c_2}{2}}\le2$, \eqref{e:onestep} follows.

Finally, the balls counted in $N_B$ are children of $B$; each is
counted at most once, because the $\A_i$ are disjoint; and $B_i$ itself is never
counted, since $B_i\subset\ov B_{\rho_i}(x_i)$ and $\rho_i\le\rb_i<r_i$. Hence
$N_B=0$ if $\#\operatorname{ch}(B)=1$, and $N_B\le\#\operatorname{ch}(B)\le
2(\#\operatorname{ch}(B)-1)$ otherwise.

We now iterate \eqref{e:onestep}. Summing \eqref{e:onestep} over $B\in\G_{k+1}$ and using~\ref{o2},
\begin{equation}\label{e:iter}
	\|\tau_{k+1}\|\le q\,\|\tau_k\|+C_1\,\mu_1(D_k)+C_1\,\eps N_k,
	\qquad D_k:=\bigcup_{B\in\G_{k+1}}D_B,\quad N_k:=\sum_{B\in\G_{k+1}}N_B .
\end{equation}
The sets $D_k$, $0\le k<K$, are pairwise disjoint and contained in $B_{3/4}$:
indeed $D_k\subset\cup \F_{k+1}$ and $D_k\cap\cup\F_k=\emptyset$, while
$\cup\F_{k+1}\subset\cup \F_j$ for $j>k$ by Proposition~\ref{p:step}(i); thus
$D_k\subset\cup \F_j$ and $D_j\cap\cup \F_j=\emptyset$. Hence
$\sum_k\mu_1(D_k)\le\mu_1(B_1)$. Moreover, since the children of the balls of
$\G_{k+1}$ are distinct elements of $\G_k$,
\[
\sum_{k=0}^{K-1}N_k\le2\sum_{k=0}^{K-1}\sum_{B\in\G_{k+1}}(\#\operatorname{ch}(B)-1)
\le2\sum_{k=0}^{K-1}\big(\#\G_k-\#\G_{k+1}\big)\le2\,\#\G_0\le2\,\#\F_0
\le\frac{2\mu_0(B_1)}{\pi\ell^2\eps}
\]
by \eqref{e:F0}. Let $T:=\sum_{k=1}^K\|\tau_k\|$. Summing \eqref{e:iter} over
$k=0,\dots,K-1$,
\[
T\le q\big(\|\tau_0\|+T\big)+C_1\mu_1(B_1)+\frac{2C_1}{\pi\ell^2}\,\mu_0(B_1),
\]
and by~\ref{o2}, \ref{o3} we conclude
\begin{equation}\label{e:tauK}
	K\,\|\tau_K\|\le T\le C_2\,\big(\mu_0(B_1)+\mu_1(B_1)\big),\qquad
	C_2:=\frac{1+C_1+2C_1/(\pi\ell^2)}{1-q}.
\end{equation}

Let $E:=\{r\in(1/2,5/8):\ \partial B_r\cap\ov B=\emptyset\ \ \forall B\in\F_K\}$.
For each $B\in\F_K$, the radii $r$ with $\partial B_r\cap\ov B\ne\emptyset$ form an
interval of length at most $2\rho_B$; hence
$|E|\ge\frac18-2\Sigma(\F_K)\ge\frac18-2\theta=\frac1{16}$. Let $r\in E$. Every
$B\in\F_K$ meeting $B_r$ is contained in $B_r\subset B_{3/4}$, hence belongs to
$\G_K$. Since $S\cap\ov B_r\subset S\cap\ov B_{7/8}\subset\cup\F_K$, we get
$\partial B_r\cap S=\emptyset$ and
$\nu(B_r)=\sum_{B\in\G_K,\,B\subset B_r}\nu(B)$. By \eqref{e:tauK} and \eqref{e:K},
\begin{align*}
|\nu(B_r)|&\le\|\tau_K\|\le\frac{C_2}{K}\big(\mu_0(B_1)+\mu_1(B_1)\big)
\le\frac{2C_2\log\Lambda}{|\log \mu_0(B_1)|}\big(\mu_0+\mu_1\big)(B_1)\\
&=2C_2\log\Lambda\;\omega(\mu_0(B_1))\,\big(\mu_0+\mu_1\big)(B_1),
\end{align*}
because $\mu_0(B_1)\le\delta_1$. Together with the trivial cases discussed at the beginning, the theorem holds with
$C:=\max\{|\log\delta_1|,\,2C_2\log\Lambda\}$. \qed

\section{Energy lower bound}\label{s:lower}

Throughout this section $\delta_1$ and $\omega$ are as in Theorem~\ref{t:log}, and we set $g(t):=t\,\omega(t)$ for $t\ge0$. The function $\omega$ is continuous and nondecreasing on $[0,\infty)$, and $g$ is a continuous increasing bijection of $[0,\infty)$ onto itself. Moreover,
\begin{equation}\label{e:g-inverse}
g^{-1}(s)\ge \tfrac12\, s\,|\log s| \qquad \text{for } s>0 \text{ small enough}.
\end{equation}
Indeed, if $t=g^{-1}(s)$ and $s$ is small, then $t<\min\{\delta_1,e^{-1}\}$, so that $s=t/|\log t|\le t$ and $t=s|\log t|\le s|\log s|\le s^{1/2}$, which gives $|\log t|\ge\frac12|\log s|$.
%

We recall the following result proven in \cite[Lemma~11]{LauteriLuckhaus2016}.

\begin{lemma}\label{l:inversion}
	Let $\Omega\subset\R^2$ be open, let $\nu$ be an $\R^m$-valued Radon measure on $\Omega$ and let $\kappa$ be a finite positive Radon measure on $\Omega$. Assume that there exists $C>0$ such that, for every $x_0\in\Omega$ and every $R>0$ with $B_R(x_0)\subset\Omega$, there exists $r\in(R/2,5R/8)$ with
	\[
	|\nu(B_r(x_0))|\le C\,\omega\Big(\frac{\kappa(B_R(x_0))}{R}\Big)\,\kappa(B_R(x_0)).
	\]
	Set
	\[
	\Theta^*(x):=\limsup_{\rho\downarrow0}\frac{\kappa(B_\rho(x))}{\rho},\qquad S:=\{x\in\Omega:\ \Theta^*(x)>0\}.
	\]
	Then:
	\begin{enumerate}
		\item $|\nu|(\Omega\setminus S)=0$;
		\item $\Hu\res S$ is $\sigma$-finite;
		\item there exists $C'>0$, depending only on $C$ and $\delta_1$, such that $|\nu|\le C'\,(\omega\circ\Theta^*)\,\kappa\res S$ in $\Omega$.
	\end{enumerate}
\end{lemma}

\subsection{Proof of Theorem~\ref{t:main}}

\begin{proof}
The upper bound follows from the construction in Appendix \ref{app:upper-bound} (see also \cite{FortunaGarroniSpadaro2025}), and we prove the lower bound. In what follows $C,c>0$ denote constants depending only on $\tau_0$ and $\lambda$, which may change from line to line; in particular they do not depend on $L$, $\delta$ and $\mu\in[\tfrac12,1)$, the uniformity in $\mu$ coming from Lemma~\ref{l:energy transfer}.

We may assume $E_L<\infty$. Let $\eps_j\downarrow0$ and $(A_j,S_j)\in\C_{\eps_j}(\mu)$ be such that $\eps_j^{-1}\mathcal{F}_{\varepsilon_j}(A_j,S_j;Q_L)\to E_L$, fix $L'\in(\delta,L)$ and let $(\widetilde A_j,\widetilde S_j)$ and $(A^{\mathrm{harm}}_j,S^{\mathrm{harm}}_j)$ be given by Lemmas~\ref{l:asymptoticequivalent} and~\ref{l:harmonic} on $Q_{L'}$ (for $j$ large $L'\le L-3\eps_j$), applied to $(A_j,S_j)\in\mathcal C_{\eps_j}(\mu)\subset\Chat_{\eps_j}$. The new pairs belong to $\Chat_{\eps_j}$, i.e.\ they satisfy \eqref{e:quantization}, but in general not \eqref{e:prescribed-phase}; the boundary condition will be used only for the original fields $A_j$, in \eqref{e:phases} below. Recall that $S^{\mathrm{harm}}_j=\widetilde S_j$ and $\curl A^{\mathrm{harm}}_j=\curl\widetilde A_j$ by \eqref{e:defAharm}.
Since $\dist_{K_\mu}(\widetilde A_j)$ is bounded by Lemma~\ref{l:asymptoticequivalent}(i), Lemma~\ref{l:asymptoticequivalent}(ii) gives
\begin{equation}\label{e:dist-tilde}
\int_{Q_{L'}}\dist^2_{K_\mu}(\widetilde A_j)\,\d x\le \mathcal{F}_{\varepsilon_j}(\widetilde A_j,\widetilde S_j;Q_{L'})+C\,\bigl|B_{\lambda\eps_j}(\widetilde S_j)\cap Q_{L'}\bigr|
\le C\,\mathcal{F}_{\varepsilon_j}(A_j,S_j;Q_L),
\end{equation}
and testing \eqref{e:pde} with its solution gives $\|A^{\mathrm{harm}}_j-\widetilde A_j\|^2_{L^2(Q_{L'})}\le C\int_{Q_{L'}}\dist^2_{K_\mu}(\widetilde A_j)$. Therefore the measures
\[
\kappa_j:=\frac1{\eps_j}\Big(\dist^2_{K_\mu}(A^{\mathrm{harm}}_j)\,\Le+\Le\res B_{\lambda\eps_j}(\widetilde S_j)\Big)\res Q_{L'}
\]
satisfy
\begin{equation}\label{e:kappa-bound}
\kappa_j(Q_{L'})\le \frac{C}{\eps_j}\mathcal{F}_{\varepsilon_j}(A_j,S_j;Q_L)
\qquad\text{and}\qquad
|\curl\widetilde A_j|\le C\,\kappa_j ,
\end{equation}
where the second inequality follows from Lemma~\ref{l:asymptoticequivalent}(iii).

We apply Theorem~\ref{t:log} in $Q_{L'}$ to $\eps_j$, to $\nu:=\frac{\lambda}{C}\curl A^{\mathrm{harm}}_j$, with $C$ the constant of Lemma~\ref{l:asymptoticequivalent}(iii), to $S:=\widetilde S_j$ and to $\sigma:=\dist^2_{K_\mu}(A^{\mathrm{harm}}_j)\,\Le\res Q_{L'}$, so that $\mu_0^{\eps_j}+\mu_1^{\eps_j}=\kappa_j$. Hypothesis (H0) holds by Lemma~\ref{l:asymptoticequivalent}(iii). As for (H1), if $\partial E$ is a closed Lipschitz curve which does not meet $\widetilde S_j$, then $A^{\mathrm{harm}}_j$ is smooth in a neighborhood of $\partial E$ and, by Stokes' theorem, $\nu(E)$ is $\pm\frac\lambda C$ times the circulation of $A^{\mathrm{harm}}_j$ along $\partial E$, so that \eqref{e:quantization}, which holds for $A^{\mathrm{harm}}_j$ since $(A^{\mathrm{harm}}_j,S^{\mathrm{harm}}_j)\in\Chat_{\eps_j}$, gives (H1) with $\tau_0$ replaced by $\lambda\tau_0/C$. Finally, if $\F$ is as in (H2), its balls cover both the support of the curl and the set where $A^{\mathrm{harm}}_j$ may fail to be harmonic, so that (H2) follows from Lemma~\ref{l:energy transfer}, after renaming $\eta$. Since the proof of Theorem~\ref{t:log} only uses the hypotheses inside $B_R(x_0)$ and $\omega$ is nondecreasing, for every $x_0\in Q_{L'}$ and every $R\ge 8\eps_j$ with $B_R(x_0)\subset Q_{L'}$ there exists a Borel set $E^j_R\subset(R/2,5R/8)$ with $|E^j_R|\ge R/16$ such that
\begin{equation}\label{e:log-j}
|\curl\widetilde A_j(B_r(x_0))|\le C\,\omega\Big(\frac{\kappa_j(B_R(x_0))}{R}\Big)\,\kappa_j(B_R(x_0))\qquad\forall\, r\in E^j_R .
\end{equation}

By \eqref{e:dist-tilde} and \eqref{e:kappa-bound}, Proposition~\ref{p:compactness} applies to $(\widetilde A_j)$: up to a subsequence, $\widetilde A_j\to A$ strongly in $L^1(Q_{L'})$, with $A\in BV(Q_{L'};K_\mu)$ satisfying \eqref{e:curlcontrolsgradient}, whose constant may depend on $\mu$ but is used only qualitatively. Up to a further subsequence, $\kappa_j\overset{*}{\rightharpoonup}\kappa$, $\curl\widetilde A_j\overset{*}{\rightharpoonup}\curl A$, and $|\curl\widetilde A_j|$ converges weakly$^*$ to a measure which is bounded by $C\kappa$ by \eqref{e:kappa-bound}; in particular $|\curl A|\le C\kappa$. By lower semicontinuity,
\begin{equation}\label{e:kappa-EL}
\kappa(Q_{L'})\le\liminf_{j\to\infty}\kappa_j(Q_{L'})\le C\,E_L .
\end{equation}
Moreover, by \eqref{e:prescribed-phase} for $A_j$ and Lemma~\ref{l:asymptoticequivalent}(ii), $\widetilde A_j-A_j\mathbbm{1}_{Q_{L'}\setminus B_{\lambda\eps_j}(S_j)}\to0$ in $L^2(Q_{L'})$, while $|B_{\lambda\eps_j}(S_j)\cap Q_L|\le C\eps_j$, so that
\begin{equation}\label{e:phases}
A=\Id\ \text{ a.e. in } Q_{L'}\cap\{x_2>\delta\},\qquad A=\mu\,\Id\ \text{ a.e. in } Q_{L'}\cap\{x_2<-\delta\}.
\end{equation}

We now pass to the limit in \eqref{e:log-j}. Fix $x_0\in Q_{L'}$ and consider the radii $\rho<\dist(x_0,\partial Q_{L'})$ with $\kappa(\partial B_\rho(x_0))=0$, which are all but countably many. For such radii $\kappa_j(B_\rho(x_0))\to\kappa(B_\rho(x_0))$ and, since the weak$^*$ limit of $|\curl\widetilde A_j|$ does not charge $\partial B_\rho(x_0)$, also $\curl\widetilde A_j(B_\rho(x_0))\to\curl A(B_\rho(x_0))$. If $R$ is such a radius, by the reverse Fatou lemma the set of radii belonging to $E^j_R$ for infinitely many $j$ has measure at least $R/16$, and for each of its elements $r$ with $\kappa(\partial B_r(x_0))=0$ we can pass to the limit in \eqref{e:log-j} along a subsequence, using the continuity of $\omega$:
\begin{equation}\label{e:log-limit}
|\curl A(B_r(x_0))|\le C\,\omega\Big(\frac{\kappa(B_R(x_0))}{R}\Big)\,\kappa(B_R(x_0)).
\end{equation}
For a general $R<\dist(x_0,\partial Q_{L'})$, we take radii $R_k\uparrow R$ as above. The radii belonging to the corresponding sets for infinitely many $k$ form a set of measure at least $R/16$ contained in $[R/2,5R/8]$, and for each such $r\in(R/2,5R/8)$ the monotonicity of $\omega$ gives $|\curl A(B_r(x_0))|\le C\omega(\kappa(B_R(x_0))/R_k)\,\kappa(B_R(x_0))$ for infinitely many $k$; letting $k\to\infty$ we obtain \eqref{e:log-limit}. Hence Lemma~\ref{l:inversion} applies to $\nu=\curl A$ and $\kappa$ in $\Omega=Q_{L'}$: with $\Theta^*$ and $S$ as there, $\Hu\res S$ is $\sigma$-finite, $|\curl A|$ is concentrated on $S$ and
\begin{equation}\label{e:curl-density}
|\curl A|\le C\,(\omega\circ\Theta^*)\,\kappa\res S .
\end{equation}

By \eqref{e:curlcontrolsgradient}, also $|DA|$ is concentrated on $S$. Since $\Le(S)=0$ and the Cantor part of the derivative of a $BV$ function vanishes on $\Hu$-$\sigma$-finite sets \cite[Proposition~3.92]{AmbrosioFuscoPallara2000}, we deduce that $A\in SBV(Q_{L'};K_\mu)$ and $DA=(A^+-A^-)\otimes n_A\,\Hu\res J_A$, so that $\curl A=(A^+-A^-)\,t_A\,\Hu\res J_A$ with $t_A$ a unit tangent to $J_A$. As $K_\mu$ is contained in the linear space $CO_+$, the jump $A^+-A^-$ is conformal and
\begin{equation}\label{e:curl-jump}
|\curl A|=\frac{|A^+-A^-|}{\sqrt2}\,\Hu\res J_A .
\end{equation}
Let $f$ be the density of the absolutely continuous part of $\kappa\res J_A$ with respect to $\Hu\res J_A$. Then $\Theta^*=2f$ at $\Hu$-a.e.\ point of $J_A$: since $J_A$ is countably $\Hu$-rectifiable, by the Besicovitch differentiation theorem the absolutely continuous part has density $2f$ and the singular part has density zero $\Hu$-a.e.\ on $J_A$, while the part of $\kappa$ outside $J_A$ has zero upper density $\Hu$-a.e.\ on $J_A$ by the density theorem \cite[Theorem~2.56]{AmbrosioFuscoPallara2000}. Comparing the densities with respect to $\Hu\res J_A$ in \eqref{e:curl-density}, by \eqref{e:curl-jump} we get $|A^+-A^-|\le C f\,\omega(2f)=\frac C2\,g(2f)$ $\Hu$-a.e.\ on $J_A$, that is, for a suitable $c>0$,
\begin{equation}\label{e:jump-energy}
\int_{J_A}g^{-1}\bigl(c|A^+-A^-|\bigr)\,\d\Hu\le 2\int_{J_A}f\,\d\Hu\le 2\,\kappa(Q_{L'}).
\end{equation}

Finally, since $|P|=\sqrt2$ for $P\in \SO(2)$ and $|P|=\sqrt2\mu$ for $P\in\mu \SO(2)$, the function
\[
\chi:=\frac{|A|/\sqrt2-\mu}{1-\mu}
\]
belongs to $BV(Q_{L'};\{0,1\})$, with $J_\chi\subset J_A$ up to an $\Hu$-null set and $|A^+-A^-|\ge\bigl||A^+|-|A^-|\bigr|=\sqrt2\,\beta$ on $J_\chi$. By \eqref{e:phases}, $\chi=1$ above $\delta$ and $\chi=0$ below $-\delta$, so that almost every vertical slice of $\chi$ has variation at least one and $\Hu(J_\chi\cap Q_{L'})\ge|D_2\chi|(Q_{L'})\ge 2L'$. Since $g^{-1}$ is increasing, by \eqref{e:jump-energy}, \eqref{e:g-inverse} and \eqref{e:kappa-EL} we conclude that, for $\beta$ small enough (which fixes $\mu_0$),
\[
c\,L'\beta|\log\beta|\le 2L'\,g^{-1}(\sqrt2\,c\,\beta)\le\int_{J_A}g^{-1}\bigl(c|A^+-A^-|\bigr)\,\d\Hu\le2\kappa(Q_{L'})\le C\,E_L .
\]
Letting $L'\uparrow L$ we obtain $E_L\ge c\,L\beta|\log\beta|$.
\end{proof}

	\begin{remark}
	The logarithmic estimate above is independent of the boundary
	data, which enter only in the final step of the proof. For
	boundary values $R_{\beta}$ and $\mu R_\alpha$, the tangential
	mismatch is
	\[
	d_{\mu,\alpha,\beta}
	:= |(R_{\beta}-\mu R_\alpha)e_1|
	= \sqrt{(1-\mu)^2+2\mu(1-\cos(\alpha-\beta))},
	\]
	and a similar argument yields
	\[
	\liminf_{\varepsilon\to0}\frac{1}{\varepsilon}
	\mathcal{F}_{\varepsilon}(A_\varepsilon,S_\varepsilon;Q_L)
	\geq cL\,d_{\mu,\alpha,\beta}(1+|\log d_{\mu,\alpha,\beta}|).
	\]
	For $\mu$ close to $1$ and $|\alpha-\beta|$ small, this is comparable to
	\[
	L\Big[(1-\mu)(1+|\log(1-\mu)|)
	+|\alpha-\beta|(1+|\log|\alpha-\beta||)\Big],
	\]
	with vanishing terms understood by continuity.
	Thus the analysis recovers the van der Merwe and Read--Shockley
	scalings, including the Lauteri--Luckhaus case $\mu=1$. This means that, up to multiplicative constants, the energetic cost
	is the same whether mismatch and misorientation are accommodated
	at a single interface or separately at two interfaces.
\end{remark}

\section*{Acknowledgements}
The authors wish to thank Marco Cicalese, Leonard Kreutz and Tim Heilmann for useful discussions. V.C. thanks the Technical University of Munich for its hospitality during part of this work. V.C. is a member of INdAM-GNAMPA and has been partially supported by the INdAM-GNAMPA project E53C25002010001; V.C. also acknowledges the support of the Sapienza grant “Avvio alla Ricerca 2025”.

\section*{Declaration of AI usage}
The authors acknowledge the use of AI tools in the preparation of this manuscript. These tools were used to revise the text and improve its exposition. The fundamental mathematical ideas and results are the work of the authors.
\appendix

\section{Blow-up lemma }

We prove here a lemma on the compactness of blow-ups of Lipschitz graph domains (see also \cite{DeLellisSzekelyhidi2006}).

We write points of $\R^n$ as $y=(y',y_n)$ with $y'\in\R^{n-1}$.

\begin{lemma}\label{l:C}
	Let $r_0>0$ and $\Omega=\{(x',x_n):x'\in(0,r_0)^{n-1},\ f(x')<x_n<r_0\}$, where
	$f:[0,r_0]^{n-1}\to\R$ is Lipschitz with $\max f<r_0$. Set $\Lambda:={\rm Lip}(f)$,
	\[
	L:=\max\{1,\,4\sqrt{n-1}\,\Lambda\},\qquad Q^L_1:=(-1,1)^{n-1}\times(-L,L),\qquad
	Q^L_r(x):=x+rQ^L_1 .
	\]
	For every $R_0>0$ there exists $\gamma=\gamma(\Omega,R_0)>0$ with the following property. Let $x_j\in\Omega$, $r_j\in(0,R_0]$ and $U_j:=r_j^{-1}(\Omega-x_j)\cap Q^L_1$. Then $|U_j|\ge\gamma$ for every $j$ and, up to a subsequence, there is a connected open set $U\subset Q^L_1$ such that $|U_j\triangle U|\to0$ and, for every $\delta>0$, there is a connected open set $U'\subset U$ with $|U\setminus U'|<\delta$ and $U'\subset U_j$ for all $j$	large enough.
\end{lemma}

\begin{proof}[Proof of Lemma \ref{l:C}]
There is no loss of generality in assuming that $r_0=1$, because the statement is invariant under homothety. 
	Write $x_j=(x_j',t_j)$. Since $x_j\in\Omega$,
	\[
	r_j^{-1}(\Omega-x_j)=\big\{y:\ y'\in r_j^{-1}\big((0,1)^{n-1}-x_j'\big),\
	f_j(y')<y_n<h_j\big\},
	\]
	where
	\[
	f_j(y'):=\frac{f(x_j'+r_jy')-t_j}{r_j},\qquad h_j:=\frac{1-t_j}{r_j}>0 .
	\]
	Therefore
	\[
	U_j=\big\{y:\ y'\in P_j,\ \varphi_j(y')<y_n<\psi_j\big\},
	\]
	with
	\[
	P_j:=r_j^{-1}\big((0,1)^{n-1}-x_j'\big)\cap(-1,1)^{n-1},\qquad
	\varphi_j:=\max\{f_j,-L\},\qquad\psi_j:=\min\{h_j,L\}.
	\]
	We collect three elementary facts.
	
	(a) $P_j=\prod_{i=1}^{n-1}(a_{j,i},b_{j,i})$ with $-1\le a_{j,i}\le0\le b_{j,i}\le1$
	and 
	\[\
	b_{j,i}-a_{j,i}\geq\min\{1,1/r_j\}\geq\rho:=\min\{1,1/R_0\}.
	\]
	
	(b) $\varphi_j$ is $\Lambda$-Lipschitz on $\overline{P_j}$ and, since $f_j(0)<0$
	and $|y'|\le\sqrt{n-1}$, we have 
	\[
	-L\le\varphi_j\le\sqrt{n-1}\,\Lambda\le L/4.
	\]
	
	(c) $\psi_j-\varphi_j\ge c_0:=\min\{3L/4,(1-\max f)/R_0\}$ on $\overline{P_j}$.
	Indeed, if $\psi_j=L$ this follows from (b); if $\psi_j=h_j<L$, then
	$h_j-\varphi_j=\min\{h_j-f_j,\,h_j+L\}$, with
	\[
	h_j-f_j(y')=\big(1-f(x_j'+r_jy')\big)/r_j\ge(1-\max f)/R_0\quad \text{and} \quad h_j+L>L.
	\]
	
	By (a) and (c), $|U_j|\ge c_0\rho^{n-1}=:\gamma$. Now, let $\bar\varphi_j:=\varphi_j\circ\pi_j$ denote the Lipschitz extension of $\varphi_j$ to $[-1,1]^{n-1}$, where $\pi_j$ is the nearest-point projection onto the convex set $\overline{P_j}$. Note that the extension $\bar\varphi_j$	is still $\Lambda$-Lipschitz with values in $[-L,L]$. By (a), the Arzel\`a--Ascoli theorem and the boundedness of $\psi_j\in(0,L]$, up to a subsequence
	\[
	a_{j,i}\to a_i,\quad b_{j,i}\to b_i,\quad \bar\varphi_j\to\bar\varphi\ \text{uniformly on }[-1,1]^{n-1},
	\quad \psi_j\to\psi\in[0,L],
	\]
	with $b_i-a_i\ge\rho$ and, passing to the limit in (c),
	$\psi-\bar\varphi\ge c_0$ on $\overline P$, where $P:=\prod_i(a_i,b_i)$. Set
	\[
	U:=\big\{y:\ y'\in P,\ \bar\varphi(y')<y_n<\psi\big\}.
	\]
	Every $y\in Q^L_1$ with $y'\notin\partial P$, $y_n\neq\bar\varphi(y')$ and
	$y_n\neq\psi$ belongs either to $U_j$ for all large $j$ or to
	$Q_1^L\setminus U_j$ for all large $j$, according to whether it belongs to $U$ or
	not. The excluded set is Lebesgue negligible, so $\chi_{U_j}\to\chi_U$ a.e.\
	and $|U_j\triangle U|\to0$ by dominated convergence.
	
	For $\varepsilon\in(0,\min\{\rho,c_0\}/3)$ let
	\[
	U'_\varepsilon:=\big\{y:\ y'\in P^\varepsilon,\ \bar\varphi(y')+\varepsilon<y_n<\psi-\varepsilon\big\},
	\qquad P^\varepsilon:=\prod_i(a_i+\varepsilon,b_i-\varepsilon).
	\]
	The sets $U$ and $U'_\varepsilon$ are open and connected: each is the region
	between two continuous graphs over a convex base, with nonempty vertical fibers
	(because $\psi-\bar\varphi\ge c_0>2\varepsilon$), and every point can be joined
	vertically to the graph of the continuous function $(\bar\varphi+\psi)/2$, which is
	connected. For $j$ large we have $a_{j,i}<a_i+\varepsilon$,
	$b_{j,i}>b_i-\varepsilon$, $\|\bar\varphi_j-\bar\varphi\|_\infty<\varepsilon$ and
	$|\psi_j-\psi|<\varepsilon$, hence $U'_\varepsilon\subset U_j$. Finally
	$|U\setminus U'_\varepsilon|\le2L\,|P\setminus P^\varepsilon|+2\varepsilon|P|
	\le C(n,L)\,\varepsilon$, which is smaller than $\delta$ for $\varepsilon$
	small, and we take $U':=U'_\varepsilon$.\qedhere
\end{proof}

\begin{remark}
	The shape of $Q_1^L$ matters: its height $L$ exceeds the oscillation $2\sqrt{n-1}\,\Lambda$ of a $\Lambda$-Lipschitz function on $(-1,1)^{n-1}$.
	In this way, every vertical fiber of $U_j$ is nonempty, thus ensuring the connectedness of the limit.
\end{remark}

%

\section{Besicovitch-type selection}

\begin{proof}[Proof of Lemma~\ref{l:besi}]
	For $i\in I$ set $d_i:=R-|x_i-y|$ and $R_i:=d_i/8$. By assumption $d_i\ge16\rho_i$,
	i.e.\ $\rho_i\le R_i/2$; moreover $3R_i<d_i$, so $B_{3R_i}(x_i)\subset B_R(y)$.
	Since $x\mapsto R-|x-y|$ is $1$-Lipschitz,
	\begin{equation}\label{e:slow}
		|R_i-R_j|\le\tfrac18|x_i-x_j|\qquad\forall\,i,j\in I .
	\end{equation}
	
	Let $J\subset I$ be maximal with the property
	$|x_j-x_l|\ge\frac25\max\{R_j,R_l\}$ for all distinct $j,l\in J$.
	
	We claim $B_{\rho_i}(x_i)\subset\bigcup_{j\in J}B_{R_j}(x_j)$
	for every $i\in I$. Let $i\in I\setminus J$; by maximality there is $j\in J$ with
	$|x_i-x_j|<\frac25\max\{R_i,R_j\}$.
	If $R_i\ge R_j$, then by \eqref{e:slow} $R_j\ge(1-\frac1{20})R_i=\frac{19}{20}R_i$,
	and $|x_i-x_j|+\rho_i<(\frac25+\frac12)R_i=\frac9{10}R_i<R_j$.
	If $R_j>R_i$, then $R_i\le\frac{21}{20}R_j$ and
	$|x_i-x_j|+\rho_i<(\frac25+\frac{21}{40})R_j=\frac{37}{40}R_j<R_j$.
	
	Next we bound the overlap. Fix $j\in J$ and let
	$L:=\{l\in J:\ B_{3R_l}(x_l)\cap B_{3R_j}(x_j)\neq\emptyset\}$. For $l\in L$,
	$|x_l-x_j|<3(R_l+R_j)$, and \eqref{e:slow} gives
	$R_l\le R_j+\frac38(R_l+R_j)$ and the symmetric inequality, i.e.\
	$\frac5{11}R_j\le R_l\le\frac{11}5R_j$. Hence $|x_l-x_j|<\frac{48}5R_j$, and for distinct
	$l,l'\in L$, $|x_l-x_{l'}|\ge\frac25\cdot\frac5{11}R_j=\frac2{11}R_j$. The balls
	$B_{R_j/11}(x_l)$, $l\in L$, are pairwise disjoint and contained in
	$B_{(48/5+1/11)R_j}(x_j)$, so that
	\[
	\#L\le\Big(11\Big(\frac{48}5+\frac1{11}\Big)\Big)^n=\Big(\frac{533}5\Big)^n<107^n .
	\]
	
	In the intersection graph of $\{B_{3R_j}(x_j)\}_{j\in J}$
	every vertex has fewer than $107^n-1$ neighbours, so a greedy colouring yields a
	partition $J=J_1\cup\dots\cup J_{107^n}$ such that, for each $c$, the balls
	$B_{3R_j}(x_j)$, $j\in J_c$, are pairwise disjoint. By the covering property above and the assumption
	on $\zeta$,
	\[
	\zeta(B_R(y))\le\sum_{j\in J}\zeta(B_{R_j}(x_j))=\sum_{c=1}^{107^n}\sum_{j\in J_c}\zeta(B_{R_j}(x_j)),
	\]
	and we choose $\tilde I:=J_c$ with $c$ maximizing the inner sum. Finally
	$R_j\ge2\rho_j$ for every $j$.
\end{proof}

\section{The upper bound construction}
\label{app:upper-bound}
\begin{lemma}\label{2d}
	There exists $C=C(\tau_0,\lambda)>0$ such that,
	for every $\mu\in(0,1)$, $\delta>0$ and $L>\delta$, there exist
	$(A_\varepsilon,S_\varepsilon)\in\mathcal C_\varepsilon(\mu)$ with identical
	square cores, alternating lower and upper cuts, and identical horizontal
	Burgers vectors, satisfying
	\[
	\limsup_{\varepsilon\downarrow0}\frac1\varepsilon
	\mathcal{F}_{\varepsilon}(A_\varepsilon,S_\varepsilon;Q_L)
	\le CL\beta(1+|\log\beta|),\qquad \beta:=1-\mu.
	\]
\end{lemma}
\begin{proof}
	We follow the construction in \cite{FortunaGarroniSpadaro2025}, adding the alternating placement of the cuts above and below the
	interface, to simulate the accommodation of both lattices.
	We omit the dependence on the parameters. Let $l>0$, $0<c<1$, and
	$Q=(-l,l)^2$, $Q'=[-cl,cl]^2$. Write
	$P_1=(l,l)$, $P_2=(l,-l)$, $P_3=(-l,-l)$, $P_4=(-l,l)$ and
	$P_1'=cP_1$, $P_2'=cP_2$, $P_1'=cP_3$, $P_1'=cP_1$.
	Set $\Delta_1=P_1P_1'P_4'P_1$, $\Delta_2=P_2P_2'P_3'P_3$,
	$\Delta_3=P_1P_2P_2'P_1'$ and $\Delta_4=P_4P_3P_3'P_4'$.
	
	Let $b^+,b^-\in\mathbb R^2$ be arbitrary translation vectors; their values
	will be chosen in the energy estimate. With $z=(x,y)$, define the map
	with a lower cut by
	\begin{equation}\label{constr}
		f(x,y)=\begin{cases}
			z,&\text{on }\Delta_1,\\
			z-\beta x e_1+b^+,&\text{on }\Delta_2\cap\{x>0\},\\
			z-\beta x e_1+b^-,&\text{on }\Delta_2\cap\{x<0\},\\
			z+\dfrac{x-y}{2x}(b^+-\beta x e_1),&\text{on }\Delta_3,\\[1mm]
			z+\dfrac{x+y}{2x}(b^--\beta x e_1),&\text{on }\Delta_4.
		\end{cases}
	\end{equation}
	The last two expressions are the linear interpolations in the $y$ direction.
	For the alternating cell, the upper translations are chosen so that the
	lateral tangential traces match those of $f$. Explicitly, set
	\begin{equation}\label{constr2}
		\widetilde f(x,y)=\begin{cases}
			z-2l\beta e_1-b^-,&\text{on }\Delta_1\cap\{x>0\},\\
			z+2l\beta e_1-b^+,&\text{on }\Delta_1\cap\{x<0\},\\
			z-\beta x e_1,&\text{on }\Delta_2,\\
			z-\dfrac\beta2(x-y)e_1-\dfrac{x+y}{2x}(2l\beta e_1+b^-),
			&\text{on }\Delta_3,\\[1mm]
			z-\dfrac\beta2(x+y)e_1+\dfrac{x-y}{2x}(2l\beta e_1-b^+),
			&\text{on }\Delta_4.
		\end{cases}
	\end{equation}
	Both maps agree across the oblique sides and have constant jumps on their
	respective cuts. Define $D=\nabla f$, $\widetilde D=\nabla\widetilde f$
	outside $Q'$ and set both fields equal to $\Id$ inside $Q'$.
	Their curl is supported in $Q'$; both fields equal $\Id$ on
	$\Delta_1$ and $\operatorname{diag}(\mu,1)$ on $\Delta_2$.
	
\input{FIG1}
	
	Differentiating \eqref{constr}--\eqref{constr2}, we obtain
	\begin{align*}
		D(l,y)e_2=\widetilde D(-l,y)e_2
		&=e_2+\frac\beta2e_1-\frac{b^+}{2l},\\
		\widetilde D(l,y)e_2=D(-l,y)e_2
		&=e_2-\frac\beta2e_1-\frac{b^-}{2l}.
	\end{align*}
	Hence consecutive cells can be glued without adding curl. Put $D_k=D$
	for even $k$ and $D_k=\widetilde D$ for odd $k$, and define the
	$4l$-periodic field and the defect set by
	\[
	A_\varepsilon(x,y)=\begin{cases}
		\Id,&y>l,\\
		D_k(x-2kl,y),&|x-2kl|<l,\ |y|<l,\\
		\mu\Id,&y<-l,
	\end{cases}
	\qquad S_\varepsilon=\bigcup_{k\in\mathbb Z}(Q'+2kl e_1).
	\]
	At $y=l$ and $y=-l$ the tangential traces are $e_1$ and $\mu e_1$,
	respectively, so these jumps are compatible as well. The elastic energy
	vanishes outside the strip. The counterclockwise circulations are
	\begin{equation}\label{eq:burgers-generic}
		\int_{\partial Q'}D\,t\,d\mathcal H^1=b^--b^+,
		\qquad
		\int_{\partial Q'}\widetilde D\,t\,d\mathcal H^1
		=-4l\beta e_1+b^+-b^-.
	\end{equation}
	Here the fields on $\partial Q'$ are understood through their exterior traces;
	these integrals also equal the respective circulations along $\partial Q$.
	
\input{FIG2}
	
	We now estimate the energy, still with generic $b^\pm$. On the lateral
	trapezoids, $|y|\le|x|$, and the derivatives of the two interpolations give
	\[
	|D-\Id|^2+|\widetilde D-\Id|^2
	\le C\beta^2+C\bigl(l^2\beta^2+|b^+|^2+|b^-|^2\bigr)
	\left(\frac1{x^2}+\frac{y^2}{x^4}\right).
	\]
	The upper and lower trapezoids contribute at most $C\beta^2l^2$.
	Integrating over $cl<|x|<l$, $|y|<|x|$, and using $\Id\in K_\mu$,
	we find, for $D_*=D$ or $\widetilde D$,
	\begin{equation}\label{computation}
		\int_{Q\setminus Q'}\operatorname{dist}^2(D_*,K_\mu)\,dx\,dy
		\le C\bigl(l^2\beta^2+|b^+|^2+|b^-|^2\bigr)(1+|\log c|).
	\end{equation}
	Now we choose
	\[
	b^+=l\beta e_1,\qquad b^-=-l\beta e_1,\qquad
	\kappa=\max\{1,\tau_0\},\qquad
	l=\frac{\kappa\varepsilon}{\beta},\qquad c=\beta.
	\]
	The lateral tangential traces then all equal $e_2$, and both cell types
	have the same horizontal circulation:
	\[
	\int_{\partial Q'}D\,t\,d\mathcal H^1
	=\int_{\partial Q'}\widetilde D\,t\,d\mathcal H^1
	=-2l\beta e_1=-2\kappa\varepsilon e_1.
	\]
	For every closed rectifiable curve $\gamma$ avoiding the cores, summing
	these circulations with their winding numbers gives
	\[
	\int_\gamma A_\varepsilon t\,d\mathcal H^1
	\in 2\kappa\varepsilon e_1\mathbb Z,
	\]
	where integrals are taken with parametrization and multiplicity.
	Every nonzero circulation therefore has modulus at least
	$2\kappa\varepsilon\ge\tau_0\varepsilon$. Also $l<\delta$ for small
	$\varepsilon$, so the pair is in $\mathcal C_\varepsilon(\mu)$.
	
	Since $|b^+|^2+|b^-|^2=2l^2\beta^2$, substitution in
	\eqref{computation} gives a cell elastic energy bounded by
	$C\varepsilon^2(1+|\log\beta|)$. Each core contributes at most
	$C(cl+\lambda\varepsilon)^2\le C\varepsilon^2$.
	To pass to $Q_L$, choose
	\[
	n_\varepsilon=\left\lceil\frac{L+l+\lambda\varepsilon}{2l}\right\rceil,
	\qquad N_\varepsilon=2n_\varepsilon+1.
	\]
	The $N_\varepsilon$ cells $Q+2kl e_1$, $|k|\le n_\varepsilon$, cover the
	part of the strip in $Q_L$ and include every core whose
	$\lambda\varepsilon$-neighborhood meets $Q_L$, since $cl<l$. Moreover,
	\[
	N_\varepsilon\le\frac Ll+\frac{\lambda\varepsilon}{l}+4,
	\qquad \varepsilon N_\varepsilon\longrightarrow\frac{L\beta}{\kappa}.
	\]
	Summing the cell estimates, including the boundary cells, gives
	\[
	\mathcal{F}_{\varepsilon}(A_\varepsilon,S_\varepsilon;Q_L)
	\le C N_\varepsilon\varepsilon^2(1+|\log\beta|).
	\]
	Division by $\varepsilon$ and passage to the limit prove the claim.
	The rescaled core contribution is $O(L\beta)$ and is included in the bound.
\end{proof}
	
	\begin{remark}
		Multiplying the competitor constructed above on the left by
		$R_\alpha$, we obtain a transition between $R_\alpha$ and
		$\mu R_\alpha$ with the same energy. Combining this construction
		with the upper bound of Lauteri--Luckhaus for a transition between
		$\Id$ and $R_\alpha$, in two adjacent strips, yields a
		competitor connecting $\Id$ and $\mu R_\alpha$. Its energy
		is bounded by the sum of the van der Merwe and Read--Shockley
		contributions, namely
		\[
		\limsup_{\varepsilon\to0}\frac{1}{\varepsilon}
		\mathcal{F}_{\varepsilon}(A_\varepsilon,S_\varepsilon;Q_L)
		\leq CL\Big[
		(1-\mu)\bigl(1+|\log(1-\mu)|\bigr)
		+|\alpha|\bigl(1+|\log|\alpha||\bigr)
		\Big],
		\]
		for sufficiently small $|\alpha|$, with the second term understood
		as zero when $\alpha=0$.
		The same upper bound can also be obtained by a single construction
		accommodating both the lattice mismatch and the misorientation,
		without introducing an intermediate region with gradient $R_\alpha$.
	\end{remark}

\bibliographystyle{plain}
\bibliography{bib}

\end{document}

%% file: comandivalerio.tex
\usepackage{indentfirst}           
\usepackage{microtype}             
\usepackage{amsmath} 
\usepackage{amssymb}
\usepackage{amsthm}
\usepackage{mathtools} 
\usepackage{esint}    
\usepackage{mathrsfs}              
\usepackage{bbold, bbm}            
\usepackage{bussproofs}            
\usepackage{verbatim}              
\usepackage{xspace}                
\usepackage{color}                 
\usepackage{wrapfig}               
\usepackage{enumitem}              
\usepackage{lettrine}              
\usepackage{hyperref}   
\usepackage{tikz}  
\usepackage{float}
\usepackage{bm}                    
\usepackage{version}               
\usepackage{algorithm}             
\usepackage[noend]{algcompatible}  
\usepackage{tikz}
\usepackage{tikz-3dplot}
\usetikzlibrary{3d,calc}
\DeclareMathOperator{\curl}{{curl}}
\DeclareMathOperator{\dist}{\textnormal{dist}}

\DeclareMathOperator{\SO}{\operatorname{SO}}

\DeclareMathOperator{\supp}{\operatorname{spt}}
\newcommand{\R}{\mathbb{R}}
\newcommand{\N}{\mathbb{N}}

\newcommand{\Y}{\mathbb{Y}}
\newcommand{\GY}{\mathbb{GY}}

%% file: FIG2.tex
\begin{figure}[ht]
		\centering
		\begin{tikzpicture}[x=.90cm,y=.70cm,>=latex,line width=.45pt]
			\fill[blue!4] (-1,-1) rectangle (11,1);
			\foreach \k in {0,...,5}{
				\pgfmathsetmacro{\xx}{2*\k}
				\pgfmathsetmacro{\sg}{(-1)^(\k+1)}
				\draw[gray] (\xx-1,-1) rectangle (\xx+1,1);
				\filldraw[fill=gray!35] (\xx-.16,-.16) rectangle (\xx+.16,.16);
				\draw[red!70!black,very thick] (\xx,\sg*.16)--(\xx,\sg);
				\ifodd\k
				\node at (\xx,-1.35) {$\widetilde D$};
				\draw[->,blue!65!black,thick] (\xx,0)--(\xx-.52,0);
				\else
				\node at (\xx,-1.35) {$D$};
				\draw[->,blue!65!black,thick] (\xx,0)--(\xx-.52,0);
				\fi
			}
			\node at (5,1.4) {$A_\varepsilon=\mathrm{Id}$};
			\node at (5,-2.05) {$A_\varepsilon=\mu\mathrm{Id}$};
			\node at (-1.4,0) {$\cdots$};\node at (11.4,0) {$\cdots$};
			\draw[<->] (-1,1.9)--(3,1.9) node[midway,above] {$4l$};
		\end{tikzpicture}
		\caption{The periodic interface: identical cores, alternating lower and upper cuts, and a period containing two cells. Arrows show the identical horizontal circulations in all cells.}
		\label{fig:periodic-interface}
	\end{figure}
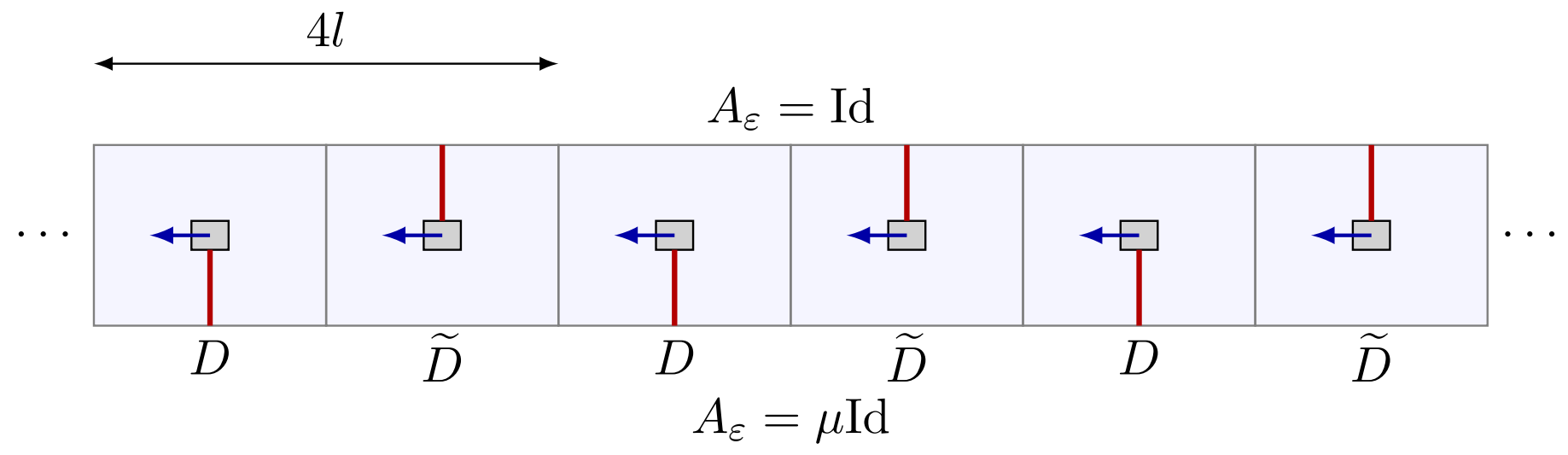